\documentclass[sn-mathphys-num]{sn-jnl}

\usepackage{amsmath,amssymb,amsfonts}%
\usepackage{amsthm}%
\usepackage{mathrsfs}%
\usepackage[title]{appendix}%
\usepackage{xcolor}%

\usepackage{tikz}
\usepackage{pgfplots}
\usetikzlibrary{calc}
\usetikzlibrary{angles,quotes} 
\usetikzlibrary{bending} 
\usetikzlibrary{arrows}
\tikzset{>=latex}

\theoremstyle{thmstyleone}%
\newtheorem{theorem}{Theorem}
\newtheorem{lemma}{Lemma}
\newtheorem{proposition}[theorem]{Proposition}%

\theoremstyle{thmstyletwo}%

\theoremstyle{thmstylethree}%
\newtheorem{definition}{Definition}%

\newcommand{\SP}[1]{\mathrm{SP}_{#1}}
\newcommand{\abs}[1]{\left\lvert #1 \right\rvert}

\newcommand{\ceil}[1]{\left\lceil #1 \right\rceil}
\newcommand{\ii}{\mathrm{i}}
\DeclareMathOperator{\cvx}{co}

\newcommand{\name}{}
\usepackage{standalone}
\usetikzlibrary{calc}
\usepackage{subcaption}

\begin{document}

\title[Zonogon approach for small polygons of maximum perimeter]%
{A zonogon approach for computing small polygons of maximum perimeter}


\author[1]{\fnm{Bernd} \sur{Mulansky}}\email{bernd.mulansky@tu-clausthal.de}
\equalcont{These authors contributed equally to this work.}
\author*[1]{\fnm{Andreas} \sur{Potschka}}\email{andreas.potschka@tu-clausthal.de}
\equalcont{These authors contributed equally to this work.}

\affil*[1]{\orgdiv{Institute of Mathematics}, \orgname{Clausthal University of Technology}, \orgaddress{\street{Erzstr.~1}, \postcode{38678} \city{Clausthal-Zellerfeld}, \country{Germany}}}


\abstract{We derive a mixed integer nonlinear programming formulation for the problem of finding a convex polygon with a given number of vertices that is small (diameter at most one) and has maximum perimeter. The formulation is based on a geometric construction using zonogons. The resulting zonogons can be characterized by an integer code and we study the number of codes that are distinct under the symmetries of the problem. We propose a two-phase computational approach. Phase~I comprises the solution of a purely combinatorial problem. Under assumption of Mossinghoff's conjecture, the Phase~I problem can be reduced to a Subset-Sum Problem. Without Mossinghoff's conjecture, a generalized Subset-Sum Problem needs to be solved, which consists of picking $n$ non-opposing $2n$-th roots of unity such that their sum is as small as possible. Phase~II consists of a non-combinatorial Nonlinear Programming Problem, which can be solved to high accuracy with arbitrary precision Newton-type methods. We provide extensive numerical results including highly accurate solutions for polygons with 64 and 128 vertices.}

\keywords{Isodiametric problem, convex polygon, maximum perimeter, sums of roots of unity}

\pacs[MSC Classification]{52-08, 52-A10, 52-A38}

\maketitle

\section{Introduction}
Let $\SP{n}$ denote the set of convex polygons in the plane with $n \ge 3$ vertices that have diameter at most one. Its elements are called \emph{small $n$-gons}. 
More than one hundred years ago, Reinhardt~\cite{reinhardt1922extremale} asked which polygons in $\SP{n}$ have maximum area or maximum perimeter. Reinhardt gave an answer to the area problem for odd $n$ (the regular $n$-gon) and for the perimeter problem for all $n$ with an odd divisor (via a construction with Reuleaux polygons). But even after a century of efforts by various bright minds, the perimeter problem is still open for the case when $n$ is a power of two.

We follow here an experimental mathematics approach, hoping that numerical solutions of a number of medium range $n$ might reveal the patterns that will eventually lead to a proof for settling Reinhardt's problems completely.
From a numerical optimization point of view, the area problem is considered an important benchmark problem, which is listed in the COPS suite~\cite{dolan2004benchmarking}.
The computational results in~\cite{graham_largest_1975} for hexagons and in~\cite{audet2002largest} for octagons, for instance, justified the conjecture that the small polygons of maximum area with an even number $n$ of vertices have a cycle of $n-1$ diameter chords with one pending diameter chord attached. This conjecture was eventually proved in~\cite{foster_diameter_2007}.

Such a result is missing for the perimeter problem, which motivates us to study novel computational methods to solve the perimeter problem for moderate sizes of $n$.

The small hexagon of largest area is well-known since its appearance in \cite{graham_largest_1975}, or, perhaps, even since \cite{bieri1961ungeloeste}. The value of its area turns out to be an algebraic number, namely a root of a polynomial of degree ten with integer coefficients. The small octagon of largest area was determined in~\cite{audet2002largest}. As established later in~\cite{audet2020using}, the maximum area of a small octagon is a root of an integer polynomial of degree $42$. 

These results suggest asking for similar representations of the longest perimeters. The longest perimeter of a small octagon was computed in~\cite{audet2002largest}, but it was already stated to four correct digits in~\cite{griffiths1975pi}. The representation as an algebraic number was still missing.
We show that the square of the longest perimeter of a small octagon is a root of an integer polynomial of degree $48$. 

For the perimeter $p(P_n)$ of a polygon $P_n \in \SP{n}$, Reinhardt~\cite{reinhardt1922extremale} derived the upper bound
\begin{equation*} 
  p(P_n) \le 2 n \sin \frac{\pi}{2n} =: \bar{u}_n < \pi.
\end{equation*}
For the case of $n$ having an odd divisor $d \ge 3$, he proposed a construction based on Reuleaux $d$-gons to show that this upper bound can be attained by polygons with a rotational $d$-symmetry. Additional global maxima were discovered in~\cite{hare_sporadic_2013} and baptized ``sporadic''. Ironically, it turned out that these non-Reinhardtian solutions are actually more frequent than Reinhardt's for larger $n$~\cite{hare2019most}.

The missing cases $n = 2^s$, $s \ge 2$, have so far only been solved for $n = 4$~\cite{taylor_simple_1953,tamvakis1987perimeter} and $n = 8$~\cite{audet2007small}.

Mossinghoff's entertaining article in the American Mathematical Monthly~\cite{mossinghoff2006problem} in 2006, together with his paper~\cite{mossinghoff2006isodiametric} helped to spark interest in the perimeter problem again. Before 2006, the interest in the problem has been rather infrequent, with a few notable contributions~\cite{taylor_simple_1953,gashkov_inequalities_1985,datta1997discrete}. Fejes T\'oth has recently added a chapter on the topic to the classic~\cite{fejes_toth_finite_2023}.

A currently highly active line of research is to improve asymptotic lower bounds for the gap between $p(P_n^\ast)$ and $\bar{u}_n$ for novel constructions of $P_n^\ast \in \SP{n}$ (see, e.g.,~\cite{mossinghoff2006isodiametric,bingane_audet_2022_tight,bingane2022tight,bingane_maximal_2023}), which have been improved to
\begin{equation}
  \bar{u}_n - p(P_n^\ast) \le \frac{\pi^9}{8 n^8} + O(n^{-10}).
\end{equation}
We shall see in the second column of Tab.~\ref{tab:perimeters_vary_n} that this bound overestimates the optimal gap considerably for $n=32, 64, 128$.

\subsection{Contributions}

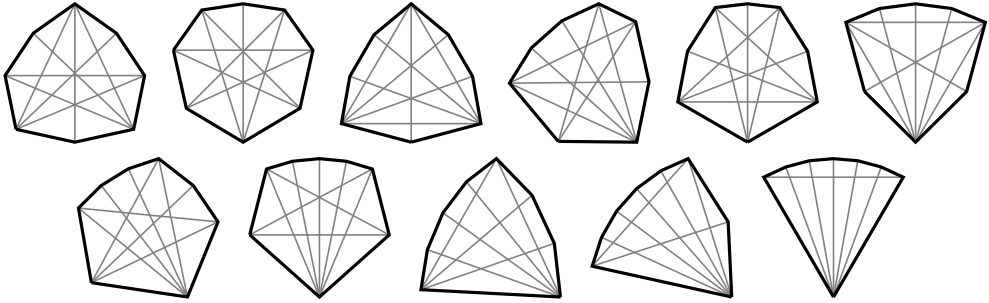
\begin{figure}[tbp]
  \centering
  \newlength{\polysize}
  \setlength\polysize{3.6cm}
\begin{tikzpicture}

\begin{axis}[
height=\polysize,
hide x axis,
hide y axis,
tick align=outside,
tick pos=left,
width=\polysize,
x grid style={white!69.0196078431373!black},
xmin=-0.549999999999438, xmax=0.550000000000562,
xtick style={color=black},
xtick={-0.6,-0.4,-0.2,0,0.2,0.4,0.6},
xticklabels={\ensuremath{-}0.6,\ensuremath{-}0.4,\ensuremath{-}0.2,0.0,0.2,0.4,0.6},
y grid style={white!69.0196078431373!black},
ymin=-0.515424160415825, ymax=0.584575839584175,
ytick style={color=black},
ytick={-0.6,-0.4,-0.2,0,0.2,0.4,0.6},
yticklabels={\ensuremath{-}0.6,\ensuremath{-}0.4,\ensuremath{-}0.2,0.0,0.2,0.4,0.6}
]
\addplot [semithick, white!50.1960784313725!black]
table {%
-0.421664714813232 -0.372175931930542
0.5 0.0158112049102783
};
\addplot [semithick, white!50.1960784313725!black]
table {%
-0.421664714813232 -0.372175931930542
0.298344135284424 0.321788907051086
};
\addplot [semithick, white!50.1960784313725!black]
table {%
-0.421664714813232 -0.372175931930542
0 0.534575819969177
};
\addplot [semithick, white!50.1960784313725!black]
table {%
-0 -0.465424180030823
0 0.534575819969177
};
\addplot [semithick, white!50.1960784313725!black]
table {%
0.421664714813232 -0.372175931930542
0 0.534575819969177
};
\addplot [semithick, white!50.1960784313725!black]
table {%
0.421664714813232 -0.372175931930542
-0.298344135284424 0.321788907051086
};
\addplot [semithick, white!50.1960784313725!black]
table {%
0.421664714813232 -0.372175931930542
-0.5 0.0158112049102783
};
\addplot [semithick, white!50.1960784313725!black]
table {%
0.5 0.0158112049102783
-0.5 0.0158112049102783
};
\addplot [very thick, black]
table {%
-0.421664714813232 -0.372175931930542
-0 -0.465424180030823
0.421664714813232 -0.372175931930542
0.5 0.0158112049102783
0.298344135284424 0.321788907051086
0 0.534575819969177
-0.298344135284424 0.321788907051086
-0.5 0.0158112049102783
-0.421664714813232 -0.372175931930542
};
\end{axis}

\end{tikzpicture}
\begin{tikzpicture}

\begin{axis}[
height=\polysize,
hide x axis,
hide y axis,
tick align=outside,
tick pos=left,
width=\polysize,
x grid style={white!69.0196078431373!black},
xmin=-0.549999999998709, xmax=0.550000000001291,
xtick style={color=black},
xtick={-0.6,-0.4,-0.2,0,0.2,0.4,0.6},
xticklabels={\ensuremath{-}0.6,\ensuremath{-}0.4,\ensuremath{-}0.2,0.0,0.2,0.4,0.6},
y grid style={white!69.0196078431373!black},
ymin=-0.641276802268483, ymax=0.458723197731517,
ytick style={color=black},
ytick={-0.8,-0.6,-0.4,-0.2,0,0.2,0.4,0.6},
yticklabels={\ensuremath{-}0.8,\ensuremath{-}0.6,\ensuremath{-}0.4,\ensuremath{-}0.2,0.0,0.2,0.4,0.6}
]
\addplot [semithick, white!50.1960784313725!black]
table {%
-0.5 0.0738769769668579
0.407464027404785 -0.346252679824829
};
\addplot [semithick, white!50.1960784313725!black]
table {%
-0.5 0.0738769769668579
0.5 0.0738769769668579
};
\addplot [semithick, white!50.1960784313725!black]
table {%
-0.407464027404785 -0.346252679824829
0.5 0.0738769769668579
};
\addplot [semithick, white!50.1960784313725!black]
table {%
-0.407464027404785 -0.346252679824829
0.296833038330078 0.363652586936951
};
\addplot [semithick, white!50.1960784313725!black]
table {%
-0 -0.591276884078979
0.296833038330078 0.363652586936951
};
\addplot [semithick, white!50.1960784313725!black]
table {%
-0 -0.591276884078979
0 0.40872323513031
};
\addplot [semithick, white!50.1960784313725!black]
table {%
-0 -0.591276884078979
-0.296833038330078 0.363652586936951
};
\addplot [semithick, white!50.1960784313725!black]
table {%
0.407464027404785 -0.346252679824829
-0.296833038330078 0.363652586936951
};
\addplot [very thick, black]
table {%
-0.5 0.0738769769668579
-0.407464027404785 -0.346252679824829
-0 -0.591276884078979
0.407464027404785 -0.346252679824829
0.5 0.0738769769668579
0.296833038330078 0.363652586936951
0 0.40872323513031
-0.296833038330078 0.363652586936951
-0.5 0.0738769769668579
};
\end{axis}

\end{tikzpicture}
\begin{tikzpicture}

\begin{axis}[
height=\polysize,
hide x axis,
hide y axis,
tick align=outside,
tick pos=left,
width=\polysize,
x grid style={white!69.0196078431373!black},
xmin=-0.55000000000159, xmax=0.54999999999841,
xtick style={color=black},
xtick={-0.6,-0.4,-0.2,0,0.2,0.4,0.6},
xticklabels={\ensuremath{-}0.6,\ensuremath{-}0.4,\ensuremath{-}0.2,0.0,0.2,0.4,0.6},
y grid style={white!69.0196078431373!black},
ymin=-0.521682885414723, ymax=0.578317114585277,
ytick style={color=black},
ytick={-0.6,-0.4,-0.2,0,0.2,0.4,0.6},
yticklabels={\ensuremath{-}0.6,\ensuremath{-}0.4,\ensuremath{-}0.2,0.0,0.2,0.4,0.6}
]
\addplot [semithick, white!50.1960784313725!black]
table {%
-0 -0.471682906150818
0 0.528317093849182
};
\addplot [semithick, white!50.1960784313725!black]
table {%
0.5 -0.337708234786987
0 0.528317093849182
};
\addplot [semithick, white!50.1960784313725!black]
table {%
0.5 -0.337708234786987
-0.26604437828064 0.305079340934753
};
\addplot [semithick, white!50.1960784313725!black]
table {%
0.5 -0.337708234786987
-0.439692616462708 0.00431180000305176
};
\addplot [semithick, white!50.1960784313725!black]
table {%
0.5 -0.337708234786987
-0.5 -0.337708234786987
};
\addplot [semithick, white!50.1960784313725!black]
table {%
0.439692616462708 0.00431180000305176
-0.5 -0.337708234786987
};
\addplot [semithick, white!50.1960784313725!black]
table {%
0.26604437828064 0.305079340934753
-0.5 -0.337708234786987
};
\addplot [semithick, white!50.1960784313725!black]
table {%
0 0.528317093849182
-0.5 -0.337708234786987
};
\addplot [very thick, black]
table {%
-0 -0.471682906150818
0.5 -0.337708234786987
0.439692616462708 0.00431180000305176
0.26604437828064 0.305079340934753
0 0.528317093849182
-0.26604437828064 0.305079340934753
-0.439692616462708 0.00431180000305176
-0.5 -0.337708234786987
-0 -0.471682906150818
};
\end{axis}

\end{tikzpicture}
\begin{tikzpicture}

\begin{axis}[
height=\polysize,
hide x axis,
hide y axis,
tick align=outside,
tick pos=left,
width=\polysize,
x grid style={white!69.0196078431373!black},
xmin=-0.592166790423953, xmax=0.507788075640656,
xtick style={color=black},
xtick={-0.6,-0.4,-0.2,0,0.2,0.4,0.6},
xticklabels={\ensuremath{-}0.6,\ensuremath{-}0.4,\ensuremath{-}0.2,0.0,0.2,0.4,0.6},
y grid style={white!69.0196078431373!black},
ymin=-0.564936003377128, ymax=0.493675887198422,
ytick style={color=black},
ytick={-0.6,-0.4,-0.2,0,0.2,0.4,0.6},
yticklabels={\ensuremath{-}0.6,\ensuremath{-}0.4,\ensuremath{-}0.2,0.0,0.2,0.4,0.6}
]
\addplot [semithick, white!50.1960784313725!black]
table {%
-0.192506313323975 -0.511337637901306
0.362067103385925 0.320797204971313
};
\addplot [semithick, white!50.1960784313725!black]
table {%
-0.192506313323975 -0.511337637901306
0.0979284048080444 0.445557117462158
};
\addplot [semithick, white!50.1960784313725!black]
table {%
0.369655132293701 -0.516817331314087
0.0979284048080444 0.445557117462158
};
\addplot [semithick, white!50.1960784313725!black]
table {%
0.369655132293701 -0.516817331314087
-0.166627049446106 0.327221393585205
};
\addplot [semithick, white!50.1960784313725!black]
table {%
0.369655132293701 -0.516817331314087
-0.386138677597046 0.137992382049561
};
\addplot [semithick, white!50.1960784313725!black]
table {%
0.369655132293701 -0.516817331314087
-0.542168855667114 -0.106235980987549
};
\addplot [semithick, white!50.1960784313725!black]
table {%
0.45779013633728 -0.097177267074585
-0.542168855667114 -0.106235980987549
};
\addplot [semithick, white!50.1960784313725!black]
table {%
0.362067103385925 0.320797204971313
-0.542168855667114 -0.106235980987549
};
\addplot [very thick, black]
table {%
-0.192506313323975 -0.511337637901306
0.369655132293701 -0.516817331314087
0.45779013633728 -0.097177267074585
0.362067103385925 0.320797204971313
0.0979284048080444 0.445557117462158
-0.166627049446106 0.327221393585205
-0.386138677597046 0.137992382049561
-0.542168855667114 -0.106235980987549
-0.192506313323975 -0.511337637901306
};
\end{axis}

\end{tikzpicture}
\begin{tikzpicture}

\begin{axis}[
height=\polysize,
hide x axis,
hide y axis,
tick align=outside,
tick pos=left,
width=\polysize,
x grid style={white!69.0196078431373!black},
xmin=-0.550000000001217, xmax=0.549999999998783,
xtick style={color=black},
xtick={-0.6,-0.4,-0.2,0,0.2,0.4,0.6},
xticklabels={\ensuremath{-}0.6,\ensuremath{-}0.4,\ensuremath{-}0.2,0.0,0.2,0.4,0.6},
y grid style={white!69.0196078431373!black},
ymin=-0.655278883653143, ymax=0.444721116346857,
ytick style={color=black},
ytick={-0.8,-0.6,-0.4,-0.2,0,0.2,0.4,0.6},
yticklabels={\ensuremath{-}0.8,\ensuremath{-}0.6,\ensuremath{-}0.4,\ensuremath{-}0.2,0.0,0.2,0.4,0.6}
]
\addplot [semithick, white!50.1960784313725!black]
table {%
-0 -0.605278968811035
0.231356263160706 0.367590188980103
};
\addplot [semithick, white!50.1960784313725!black]
table {%
-0 -0.605278968811035
0 0.394721150398254
};
\addplot [semithick, white!50.1960784313725!black]
table {%
-0 -0.605278968811035
-0.231356263160706 0.367590188980103
};
\addplot [semithick, white!50.1960784313725!black]
table {%
0.5 -0.31440544128418
-0.231356263160706 0.367590188980103
};
\addplot [semithick, white!50.1960784313725!black]
table {%
0.5 -0.31440544128418
-0.430418252944946 0.0520941019058228
};
\addplot [semithick, white!50.1960784313725!black]
table {%
0.5 -0.31440544128418
-0.5 -0.31440544128418
};
\addplot [semithick, white!50.1960784313725!black]
table {%
0.430418252944946 0.0520941019058228
-0.5 -0.31440544128418
};
\addplot [semithick, white!50.1960784313725!black]
table {%
0.231356263160706 0.367590188980103
-0.5 -0.31440544128418
};
\addplot [very thick, black]
table {%
-0 -0.605278968811035
0.5 -0.31440544128418
0.430418252944946 0.0520941019058228
0.231356263160706 0.367590188980103
0 0.394721150398254
-0.231356263160706 0.367590188980103
-0.430418252944946 0.0520941019058228
-0.5 -0.31440544128418
-0 -0.605278968811035
};
\end{axis}

\end{tikzpicture}
\begin{tikzpicture}

\begin{axis}[
height=\polysize,
hide x axis,
hide y axis,
tick align=outside,
tick pos=left,
width=\polysize,
x grid style={white!69.0196078431373!black},
xmin=-0.549999999999332, xmax=0.550000000000668,
xtick style={color=black},
xtick={-0.6,-0.4,-0.2,0,0.2,0.4,0.6},
xticklabels={\ensuremath{-}0.6,\ensuremath{-}0.4,\ensuremath{-}0.2,0.0,0.2,0.4,0.6},
y grid style={white!69.0196078431373!black},
ymin=-0.724494158464487, ymax=0.375505841535514,
ytick style={color=black},
ytick={-0.8,-0.6,-0.4,-0.2,0,0.2,0.4},
yticklabels={\ensuremath{-}0.8,\ensuremath{-}0.6,\ensuremath{-}0.4,\ensuremath{-}0.2,0.0,0.2,0.4}
]
\addplot [semithick, white!50.1960784313725!black]
table {%
-0.366025447845459 -0.308468818664551
0.5 0.191531181335449
};
\addplot [semithick, white!50.1960784313725!black]
table {%
-0 -0.67449414730072
0.5 0.191531181335449
};
\addplot [semithick, white!50.1960784313725!black]
table {%
-0 -0.67449414730072
0.258819103240967 0.291431665420532
};
\addplot [semithick, white!50.1960784313725!black]
table {%
-0 -0.67449414730072
0 0.32550585269928
};
\addplot [semithick, white!50.1960784313725!black]
table {%
-0 -0.67449414730072
-0.258819103240967 0.291431665420532
};
\addplot [semithick, white!50.1960784313725!black]
table {%
-0 -0.67449414730072
-0.5 0.191531181335449
};
\addplot [semithick, white!50.1960784313725!black]
table {%
0.366025447845459 -0.308468818664551
-0.5 0.191531181335449
};
\addplot [semithick, white!50.1960784313725!black]
table {%
0.5 0.191531181335449
-0.5 0.191531181335449
};
\addplot [very thick, black]
table {%
-0.366025447845459 -0.308468818664551
-0 -0.67449414730072
0.366025447845459 -0.308468818664551
0.5 0.191531181335449
0.258819103240967 0.291431665420532
0 0.32550585269928
-0.258819103240967 0.291431665420532
-0.5 0.191531181335449
-0.366025447845459 -0.308468818664551
};
\end{axis}

\end{tikzpicture}
\begin{tikzpicture}

\begin{axis}[
height=\polysize,
hide x axis,
hide y axis,
tick align=outside,
tick pos=left,
width=\polysize,
x grid style={white!69.0196078431373!black},
xmin=-0.521561513759966, xmax=0.573148884685125,
xtick style={color=black},
xtick={-0.6,-0.4,-0.2,0,0.2,0.4,0.6},
xticklabels={\ensuremath{-}0.6,\ensuremath{-}0.4,\ensuremath{-}0.2,0.0,0.2,0.4,0.6},
y grid style={white!69.0196078431373!black},
ymin=-0.639457350479784, ymax=0.437135402013057,
ytick style={color=black},
ytick={-0.8,-0.6,-0.4,-0.2,0,0.2,0.4,0.6},
yticklabels={\ensuremath{-}0.8,\ensuremath{-}0.6,\ensuremath{-}0.4,\ensuremath{-}0.2,0.0,0.2,0.4,0.6}
]
\addplot [semithick, white!50.1960784313725!black]
table {%
-0.37955641746521 -0.489131569862366
0.523389339447021 -0.0593769550323486
};
\addplot [semithick, white!50.1960784313725!black]
table {%
-0.37955641746521 -0.489131569862366
0.347173929214478 0.197791337966919
};
\addplot [semithick, white!50.1960784313725!black]
table {%
-0.37955641746521 -0.489131569862366
0.100329399108887 0.388199329376221
};
\addplot [semithick, white!50.1960784313725!black]
table {%
0.305526494979858 -0.590521335601807
0.100329399108887 0.388199329376221
};
\addplot [semithick, white!50.1960784313725!black]
table {%
0.305526494979858 -0.590521335601807
-0.11550235748291 0.316525936126709
};
\addplot [semithick, white!50.1960784313725!black]
table {%
0.305526494979858 -0.590521335601807
-0.309558391571045 0.19793963432312
};
\addplot [semithick, white!50.1960784313725!black]
table {%
0.305526494979858 -0.590521335601807
-0.471801996231079 0.0385736227035522
};
\addplot [semithick, white!50.1960784313725!black]
table {%
0.523389339447021 -0.0593769550323486
-0.471801996231079 0.0385736227035522
};
\addplot [very thick, black]
table {%
-0.37955641746521 -0.489131569862366
0.305526494979858 -0.590521335601807
0.523389339447021 -0.0593769550323486
0.347173929214478 0.197791337966919
0.100329399108887 0.388199329376221
-0.11550235748291 0.316525936126709
-0.309558391571045 0.19793963432312
-0.471801996231079 0.0385736227035522
-0.37955641746521 -0.489131569862366
};
\end{axis}

\end{tikzpicture}
\begin{tikzpicture}

\begin{axis}[
height=\polysize,
hide x axis,
hide y axis,
tick align=outside,
tick pos=left,
width=\polysize,
x grid style={white!69.0196078431373!black},
xmin=-0.55000000000105, xmax=0.54999999999895,
xtick style={color=black},
xtick={-0.6,-0.4,-0.2,0,0.2,0.4,0.6},
xticklabels={\ensuremath{-}0.6,\ensuremath{-}0.4,\ensuremath{-}0.2,0.0,0.2,0.4,0.6},
y grid style={white!69.0196078431373!black},
ymin=-0.763853918995242, ymax=0.336146081004758,
ytick style={color=black},
ytick={-0.8,-0.6,-0.4,-0.2,0,0.2,0.4},
yticklabels={\ensuremath{-}0.8,\ensuremath{-}0.6,\ensuremath{-}0.4,\ensuremath{-}0.2,0.0,0.2,0.4}
]
\addplot [semithick, white!50.1960784313725!black]
table {%
-0.5 -0.264854907989502
0.5 -0.264854907989502
};
\addplot [semithick, white!50.1960784313725!black]
table {%
-0.5 -0.264854907989502
0.379295587539673 0.211421608924866
};
\addplot [semithick, white!50.1960784313725!black]
table {%
-0 -0.71385383605957
0.379295587539673 0.211421608924866
};
\addplot [semithick, white!50.1960784313725!black]
table {%
-0 -0.71385383605957
0.19329309463501 0.267287135124207
};
\addplot [semithick, white!50.1960784313725!black]
table {%
-0 -0.71385383605957
0 0.28614604473114
};
\addplot [semithick, white!50.1960784313725!black]
table {%
-0 -0.71385383605957
-0.19329309463501 0.267287135124207
};
\addplot [semithick, white!50.1960784313725!black]
table {%
-0 -0.71385383605957
-0.379295587539673 0.211421608924866
};
\addplot [semithick, white!50.1960784313725!black]
table {%
0.5 -0.264854907989502
-0.379295587539673 0.211421608924866
};
\addplot [very thick, black]
table {%
-0.5 -0.264854907989502
-0 -0.71385383605957
0.5 -0.264854907989502
0.379295587539673 0.211421608924866
0.19329309463501 0.267287135124207
0 0.28614604473114
-0.19329309463501 0.267287135124207
-0.379295587539673 0.211421608924866
-0.5 -0.264854907989502
};
\end{axis}

\end{tikzpicture}
\begin{tikzpicture}

\begin{axis}[
height=\polysize,
hide x axis,
hide y axis,
tick align=outside,
tick pos=left,
width=\polysize,
x grid style={white!69.0196078431373!black},
xmin=-0.527268044988985, xmax=0.571515144562132,
xtick style={color=black},
xtick={-0.6,-0.4,-0.2,0,0.2,0.4,0.6},
xticklabels={\ensuremath{-}0.6,\ensuremath{-}0.4,\ensuremath{-}0.2,0.0,0.2,0.4,0.6},
y grid style={white!69.0196078431373!black},
ymin=-0.484035894787056, ymax=0.493400889831579,
ytick style={color=black},
ytick={-0.6,-0.4,-0.2,0,0.2,0.4,0.6},
yticklabels={\ensuremath{-}0.6,\ensuremath{-}0.4,\ensuremath{-}0.2,0.0,0.2,0.4,0.6}
]
\addplot [semithick, white!50.1960784313725!black]
table {%
0.521570444107056 -0.439606904983521
0.0628466606140137 0.44897198677063
};
\addplot [semithick, white!50.1960784313725!black]
table {%
0.521570444107056 -0.439606904983521
-0.151503920555115 0.299967885017395
};
\addplot [semithick, white!50.1960784313725!black]
table {%
0.521570444107056 -0.439606904983521
-0.319985508918762 0.100563049316406
};
\addplot [semithick, white!50.1960784313725!black]
table {%
0.521570444107056 -0.439606904983521
-0.431116461753845 -0.135653495788574
};
\addplot [semithick, white!50.1960784313725!black]
table {%
0.521570444107056 -0.439606904983521
-0.477323293685913 -0.392583966255188
};
\addplot [semithick, white!50.1960784313725!black]
table {%
0.477412581443787 -0.0951292514801025
-0.477323293685913 -0.392583966255188
};
\addplot [semithick, white!50.1960784313725!black]
table {%
0.318099498748779 0.213470935821533
-0.477323293685913 -0.392583966255188
};
\addplot [semithick, white!50.1960784313725!black]
table {%
0.0628466606140137 0.44897198677063
-0.477323293685913 -0.392583966255188
};
\addplot [very thick, black]
table {%
0.521570444107056 -0.439606904983521
0.477412581443787 -0.0951292514801025
0.318099498748779 0.213470935821533
0.0628466606140137 0.44897198677063
-0.151503920555115 0.299967885017395
-0.319985508918762 0.100563049316406
-0.431116461753845 -0.135653495788574
-0.477323293685913 -0.392583966255188
0.521570444107056 -0.439606904983521
};
\end{axis}

\end{tikzpicture}
\begin{tikzpicture}

\begin{axis}[
height=\polysize,
hide x axis,
hide y axis,
tick align=outside,
tick pos=left,
width=\polysize,
x grid style={white!69.0196078431373!black},
xmin=-0.502609330507484, xmax=0.572403395735409,
xtick style={color=black},
xtick={-0.6,-0.4,-0.2,0,0.2,0.4,0.6},
xticklabels={\ensuremath{-}0.6,\ensuremath{-}0.4,\ensuremath{-}0.2,0.0,0.2,0.4,0.6},
y grid style={white!69.0196078431373!black},
ymin=-0.58472217842232, ymax=0.462829051670588,
ytick style={color=black},
ytick={-0.6,-0.4,-0.2,0,0.2,0.4,0.6},
yticklabels={\ensuremath{-}0.6,\ensuremath{-}0.4,\ensuremath{-}0.2,0.0,0.2,0.4,0.6}
]
\addplot [semithick, white!50.1960784313725!black]
table {%
0.523539185523987 -0.537106275558472
0.218436002731323 0.415213108062744
};
\addplot [semithick, white!50.1960784313725!black]
table {%
0.523539185523987 -0.537106275558472
0.0271049737930298 0.330968141555786
};
\addplot [semithick, white!50.1960784313725!black]
table {%
0.523539185523987 -0.537106275558472
-0.142529606819153 0.208784103393555
};
\addplot [semithick, white!50.1960784313725!black]
table {%
0.523539185523987 -0.537106275558472
-0.283053755760193 0.0540010929107666
};
\addplot [semithick, white!50.1960784313725!black]
table {%
0.523539185523987 -0.537106275558472
-0.388325929641724 -0.12661612033844
};
\addplot [semithick, white!50.1960784313725!black]
table {%
0.523539185523987 -0.537106275558472
-0.453745126724243 -0.325173616409302
};
\addplot [semithick, white!50.1960784313725!black]
table {%
0.498574256896973 -0.0200704336166382
-0.453745126724243 -0.325173616409302
};
\addplot [semithick, white!50.1960784313725!black]
table {%
0.218436002731323 0.415213108062744
-0.453745126724243 -0.325173616409302
};
\addplot [very thick, black]
table {%
0.523539185523987 -0.537106275558472
0.498574256896973 -0.0200704336166382
0.218436002731323 0.415213108062744
0.0271049737930298 0.330968141555786
-0.142529606819153 0.208784103393555
-0.283053755760193 0.0540010929107666
-0.388325929641724 -0.12661612033844
-0.453745126724243 -0.325173616409302
0.523539185523987 -0.537106275558472
};
\end{axis}

\end{tikzpicture}
\begin{tikzpicture}

\begin{axis}[
height=\polysize,
hide x axis,
hide y axis,
tick align=outside,
tick pos=left,
width=\polysize,
x grid style={white!69.0196078431373!black},
xmin=-0.54999999999997, xmax=0.55000000000003,
xtick style={color=black},
xtick={-0.6,-0.4,-0.2,0,0.2,0.4,0.6},
xticklabels={\ensuremath{-}0.6,\ensuremath{-}0.4,\ensuremath{-}0.2,0.0,0.2,0.4,0.6},
y grid style={white!69.0196078431373!black},
ymin=-0.872631444395639, ymax=0.227368555604361,
ytick style={color=black},
ytick={-1,-0.8,-0.6,-0.4,-0.2,0,0.2,0.4},
yticklabels={\ensuremath{-}1.0,\ensuremath{-}0.8,\ensuremath{-}0.6,\ensuremath{-}0.4,\ensuremath{-}0.2,0.0,0.2,0.4}
]
\addplot [semithick, white!50.1960784313725!black]
table {%
-0 -0.822631359100342
0.5 0.0433939695358276
};
\addplot [semithick, white!50.1960784313725!black]
table {%
-0 -0.822631359100342
0.342020153999329 0.117061138153076
};
\addplot [semithick, white!50.1960784313725!black]
table {%
-0 -0.822631359100342
0.173648118972778 0.162176370620728
};
\addplot [semithick, white!50.1960784313725!black]
table {%
-0 -0.822631359100342
0 0.177368521690369
};
\addplot [semithick, white!50.1960784313725!black]
table {%
-0 -0.822631359100342
-0.173648118972778 0.162176370620728
};
\addplot [semithick, white!50.1960784313725!black]
table {%
-0 -0.822631359100342
-0.342020153999329 0.117061138153076
};
\addplot [semithick, white!50.1960784313725!black]
table {%
-0 -0.822631359100342
-0.5 0.0433939695358276
};
\addplot [semithick, white!50.1960784313725!black]
table {%
0.5 0.0433939695358276
-0.5 0.0433939695358276
};
\addplot [very thick, black]
table {%
-0 -0.822631359100342
0.5 0.0433939695358276
0.342020153999329 0.117061138153076
0.173648118972778 0.162176370620728
0 0.177368521690369
-0.173648118972778 0.162176370620728
-0.342020153999329 0.117061138153076
-0.5 0.0433939695358276
-0 -0.822631359100342
};
\end{axis}

\end{tikzpicture}
  \caption{Eleven small octagons with locally maximal perimeter, ordered with decreasing perimeter from left to right and top to bottom. The diameter graph, also known as the skeleton of the polygon, is depicted in gray.}
  \label{fig:polys_n_8}
\end{figure}

The paper at hand will also not give a complete solution of the perimeter problem. Instead, we present a Mixed-Integer Nonlinear Programming Problem (MINLP) that is equivalent to the perimeter problem under the conjecture that optimal polygons have $n$ distinct diameters. We count the number of \emph{codes} (the integer part of the solutions) that are distinct under dihedral symmetry (see Tab.~\ref{tab:number_of_codes}). For a fixed code, the MINLP reduces to a continuous Nonlinear Programming Problem (NLP) with $n-1$ variables, two equality constraints and $n$ inequality constraints. We were not able to discover multiple local extrema for a fixed code NLP numerically, but we are not able to provide a proof for uniqueness of the NLP solution. For $n=8$, we display the local maxima for each of the eleven distinct codes in Fig.~\ref{fig:polys_n_8} along with their perimeters in Tab.~\ref{tab:perimeters_n_8}. We also enumerated all codes for the cases $n=16$ and $n=32$ and computed local maxima for each code. Fig.~\ref{fig:polys_n_32} shows the three local solutions for $n=32$ with the largest perimeters, which coincide in the first thirteen decimal places. Hence, high-precision arithmetic is required for reliably computing solutions for larger $n$.

We further propose a two-phase algorithm to compute excellent suboptimal solutions of the MINLP by solving first a linear combinatorial Subset-Sum Problem (SSP) to generate codes and, second, solving a continuous nonlinear problem with fixed code using Newton-type methods with arbitrary precision arithmetic. The derivation of the SSP exploits Mossinghoff's conjecture~\cite{mossinghoff2006isodiametric}, which says that optimal small $n$-gons are axially symmetric. The resulting solutions for $n \le 32$ coincide with the best local maxima in the full enumeration approach, and we confirm the current record of~\cite{bingane_maximal_2023}, but are able to compute the perimeter to arbitrary precision for the first time. We also provide a candidate polygon with large perimeter, which is as close as $2 \cdot 10^{-38}$ to the upper bound $\bar{u}_{128}$. Such computations were not possible with the computational methods available before. The resulting polygons are depicted in Fig.~\ref{fig:polys_vary_n} along with the distance of their respective perimeters to the upper bound $\bar{u}_n$ in Tab.~\ref{tab:perimeters_vary_n}.

Moreover, we show that the maximum perimeter of the small octagon is an algebraic number.

\begin{table}[tbp]
  \centering
  \begin{tabular}{rl}
    $n$ & \#Codes \\
    \hline
    $4$ & $1$ \\
    $8$ & $11$ \\
    $16$ & $1,087$ \\
    $32$ & $33,570,815 \approx 3.4 \cdot 10^{7}$ \\
    $64$ & $72,057,595,111,669,759 \approx 7.2 \cdot 10^{16}$ \\
    $128$ & $664,613,997,892,457,941,063,589,548,567,560,191 \approx 6.6 \cdot 10^{35}$ \\
    \hline
  \end{tabular}
  \caption{The number of distinct codes for the MINLPs~\eqref{eqn:minlp_general} and~\eqref{eqn:minlp} for $n=2^s$.}
  \label{tab:number_of_codes}
\end{table}

\begin{table}[htbp]
  \centering
  \begin{tabular}{ccc}
    Rank & Perimeter \\ \hline
    1  &  & 3.121147134059831353864659503638086530909542\dots \\ 
    2  & & 3.119597665200247590150972423994095000480919\dots \\ 
    3  & $12 \sin \frac{\pi}{18} + 4 \sin \frac{\pi}{12}$ & 3.119054312413247235616194871727970865400783\dots \\ 
    4  & & 3.116482146091382523455235401221637774453205\dots \\ 
    5  & & 3.114973336127984895463908314651370982428416\dots \\ 
    6  & $8 \sin \frac{\pi}{24} + 8 \sin \frac{\pi}{12}$ & 3.114761898580578831178440524156298708342476\dots \\ 
    7  & & 3.108103162518355196717979437084906769281062\dots \\ 
    8  & & 3.103535201958031713403480443438109805364658\dots \\ 
    9  & $1 + 6 \sin \frac{\pi}{18} + 8 \sin \frac{\pi}{24}$ & 3.086098603761994825497549583779800857552179\dots \\ 
    10 & $1 + 4 \sin \frac{\pi}{12} + 10 \sin \frac{\pi}{30}$ & 3.080560813086617763393936898521174504202834\dots \\ 
    11 & $2 + 12 \sin \frac{\pi}{36}$ & 3.045868912971898082696771250049682616532413\dots \\ 
  \end{tabular}
  \caption{Perimeters of the small polygons in Fig.~\ref{fig:polys_n_8} with exact formulas for the cases where the octagon can be constructed from an equilateral triangle.}
  \label{tab:perimeters_n_8}
\end{table}

\begin{figure}[htbp]
  \centering  
  \setlength\polysize{5.6cm}
  \begin{tabular}{ccc}
\begin{tikzpicture}

\definecolor{darkgray176}{RGB}{176,176,176}
\definecolor{gray}{RGB}{128,128,128}

\begin{axis}[
height=\polysize,
hide x axis,
hide y axis,
tick align=outside,
tick pos=left,
width=\polysize,
x grid style={darkgray176},
xmin=-0.55, xmax=0.55,
xtick style={color=black},
y grid style={darkgray176},
ymin=-0.05, ymax=1.05,
ytick style={color=black}
]
\addplot [semithick, gray]
table {%
0 0
0 1
};
\addplot [semithick, gray]
table {%
0.0980162620544434 0.00481522083282471
0 1
};
\addplot [semithick, gray]
table {%
0.0980162620544434 0.00481522083282471
-0.0970739126205444 0.985600471496582
};
\addplot [semithick, gray]
table {%
0.193209886550903 0.0286599397659302
-0.0970739126205444 0.985600471496582
};
\addplot [semithick, gray]
table {%
0.285607814788818 0.0617202520370483
-0.0970739126205444 0.985600471496582
};
\addplot [semithick, gray]
table {%
0.285607814788818 0.0617202520370483
-0.185787916183472 0.943642139434814
};
\addplot [semithick, gray]
table {%
0.285607814788818 0.0617202520370483
-0.269962072372437 0.893190145492554
};
\addplot [semithick, gray]
table {%
0.36443042755127 0.120178937911987
-0.269962072372437 0.893190145492554
};
\addplot [semithick, gray]
table {%
0.36443042755127 0.120178937911987
-0.342676043510437 0.827286124229431
};
\addplot [semithick, gray]
table {%
0.430333852767944 0.192892074584961
-0.342676043510437 0.827286124229431
};
\addplot [semithick, gray]
table {%
0.430333852767944 0.192892074584961
-0.401135444641113 0.748462677001953
};
\addplot [semithick, gray]
table {%
0.480785250663757 0.27706503868103
-0.401135444641113 0.748462677001953
};
\addplot [semithick, gray]
table {%
0.480785250663757 0.27706503868103
-0.44309401512146 0.659749150276184
};
\addplot [semithick, gray]
table {%
0.480785250663757 0.27706503868103
-0.476154923439026 0.567350149154663
};
\addplot [semithick, gray]
table {%
0.480785250663757 0.27706503868103
-0.5 0.472155451774597
};
\addplot [semithick, gray]
table {%
0.495184659957886 0.374138355255127
-0.5 0.472155451774597
};
\addplot [semithick, gray]
table {%
0.5 0.472155451774597
-0.5 0.472155451774597
};
\addplot [semithick, gray]
table {%
0.5 0.472155451774597
-0.495184659957886 0.374138355255127
};
\addplot [semithick, gray]
table {%
0.5 0.472155451774597
-0.480785250663757 0.27706503868103
};
\addplot [semithick, gray]
table {%
0.476154923439026 0.567350149154663
-0.480785250663757 0.27706503868103
};
\addplot [semithick, gray]
table {%
0.44309401512146 0.659749150276184
-0.480785250663757 0.27706503868103
};
\addplot [semithick, gray]
table {%
0.401135444641113 0.748462677001953
-0.480785250663757 0.27706503868103
};
\addplot [semithick, gray]
table {%
0.401135444641113 0.748462677001953
-0.430333852767944 0.192892074584961
};
\addplot [semithick, gray]
table {%
0.342676043510437 0.827286124229431
-0.430333852767944 0.192892074584961
};
\addplot [semithick, gray]
table {%
0.342676043510437 0.827286124229431
-0.36443042755127 0.120178937911987
};
\addplot [semithick, gray]
table {%
0.269962072372437 0.893190145492554
-0.36443042755127 0.120178937911987
};
\addplot [semithick, gray]
table {%
0.269962072372437 0.893190145492554
-0.285607814788818 0.0617202520370483
};
\addplot [semithick, gray]
table {%
0.185787916183472 0.943642139434814
-0.285607814788818 0.0617202520370483
};
\addplot [semithick, gray]
table {%
0.0970739126205444 0.985600471496582
-0.285607814788818 0.0617202520370483
};
\addplot [semithick, gray]
table {%
0.0970739126205444 0.985600471496582
-0.193209886550903 0.0286599397659302
};
\addplot [semithick, gray]
table {%
0.0970739126205444 0.985600471496582
-0.0980162620544434 0.00481522083282471
};
\addplot [semithick, gray]
table {%
0 1
-0.0980162620544434 0.00481522083282471
};
\addplot [very thick, black]
table {%
0 0
0.0980162620544434 0.00481522083282471
0.193209886550903 0.0286599397659302
0.285607814788818 0.0617202520370483
0.36443042755127 0.120178937911987
0.430333852767944 0.192892074584961
0.480785250663757 0.27706503868103
0.495184659957886 0.374138355255127
0.5 0.472155451774597
0.476154923439026 0.567350149154663
0.44309401512146 0.659749150276184
0.401135444641113 0.748462677001953
0.342676043510437 0.827286124229431
0.269962072372437 0.893190145492554
0.185787916183472 0.943642139434814
0.0970739126205444 0.985600471496582
0 1
-0.0970739126205444 0.985600471496582
-0.185787916183472 0.943642139434814
-0.269962072372437 0.893190145492554
-0.342676043510437 0.827286124229431
-0.401135444641113 0.748462677001953
-0.44309401512146 0.659749150276184
-0.476154923439026 0.567350149154663
-0.5 0.472155451774597
-0.495184659957886 0.374138355255127
-0.480785250663757 0.27706503868103
-0.430333852767944 0.192892074584961
-0.36443042755127 0.120178937911987
-0.285607814788818 0.0617202520370483
-0.193209886550903 0.0286599397659302
-0.0980162620544434 0.00481522083282471
0 0
};
\end{axis}

\end{tikzpicture} &
\begin{tikzpicture}

\begin{axis}[
height=\polysize,
hide x axis,
hide y axis,
tick align=outside,
tick pos=left,
width=\polysize,
x grid style={white!69.0196078431373!black},
xmin=-0.55000000002868, xmax=0.549999999971321,
xtick style={color=black},
xtick={-0.6,-0.4,-0.2,0,0.2,0.4,0.6},
xticklabels={\ensuremath{-}0.6,\ensuremath{-}0.4,\ensuremath{-}0.2,0.0,0.2,0.4,0.6},
y grid style={white!69.0196078431373!black},
ymin=-0.543979749575227, ymax=0.556020250424774,
ytick style={color=black},
ytick={-0.6,-0.4,-0.2,0,0.2,0.4,0.6},
yticklabels={\ensuremath{-}0.6,\ensuremath{-}0.4,\ensuremath{-}0.2,0.0,0.2,0.4,0.6}
]
\addplot [semithick, white!50.1960784313725!black]
table {%
-0.476155042648315 -0.198433399200439
0.480785131454468 0.0918514728546143
};
\addplot [semithick, white!50.1960784313725!black]
table {%
-0.476155042648315 -0.198433399200439
0.447724223136902 0.184250593185425
};
\addplot [semithick, white!50.1960784313725!black]
table {%
-0.476155042648315 -0.198433399200439
0.405765771865845 0.272964239120483
};
\addplot [semithick, white!50.1960784313725!black]
table {%
-0.425703644752502 -0.282606363296509
0.405765771865845 0.272964239120483
};
\addplot [semithick, white!50.1960784313725!black]
table {%
-0.425703644752502 -0.282606363296509
0.347306132316589 0.351787686347961
};
\addplot [semithick, white!50.1960784313725!black]
table {%
-0.359800338745117 -0.355319380760193
0.347306132316589 0.351787686347961
};
\addplot [semithick, white!50.1960784313725!black]
table {%
-0.359800338745117 -0.355319380760193
0.274592041969299 0.417691707611084
};
\addplot [semithick, white!50.1960784313725!black]
table {%
-0.280977964401245 -0.413777947425842
0.274592041969299 0.417691707611084
};
\addplot [semithick, white!50.1960784313725!black]
table {%
-0.280977964401245 -0.413777947425842
0.190417766571045 0.468143820762634
};
\addplot [semithick, white!50.1960784313725!black]
table {%
-0.192265510559082 -0.455735802650452
0.190417766571045 0.468143820762634
};
\addplot [semithick, white!50.1960784313725!black]
table {%
-0.192265510559082 -0.455735802650452
0.0980181694030762 0.50120484828949
};
\addplot [semithick, white!50.1960784313725!black]
table {%
-0.0970721244812012 -0.479580402374268
0.0980181694030762 0.50120484828949
};
\addplot [semithick, white!50.1960784313725!black]
table {%
-0 -0.493979692459106
0.0980181694030762 0.50120484828949
};
\addplot [semithick, white!50.1960784313725!black]
table {%
-0 -0.493979692459106
0 0.506020307540894
};
\addplot [semithick, white!50.1960784313725!black]
table {%
-0 -0.493979692459106
-0.0980181694030762 0.50120484828949
};
\addplot [semithick, white!50.1960784313725!black]
table {%
0.0970721244812012 -0.479580402374268
-0.0980181694030762 0.50120484828949
};
\addplot [semithick, white!50.1960784313725!black]
table {%
0.192265510559082 -0.455735802650452
-0.0980181694030762 0.50120484828949
};
\addplot [semithick, white!50.1960784313725!black]
table {%
0.192265510559082 -0.455735802650452
-0.190417766571045 0.468143820762634
};
\addplot [semithick, white!50.1960784313725!black]
table {%
0.280977964401245 -0.413777947425842
-0.190417766571045 0.468143820762634
};
\addplot [semithick, white!50.1960784313725!black]
table {%
0.280977964401245 -0.413777947425842
-0.274592041969299 0.417691707611084
};
\addplot [semithick, white!50.1960784313725!black]
table {%
0.359800338745117 -0.355319380760193
-0.274592041969299 0.417691707611084
};
\addplot [semithick, white!50.1960784313725!black]
table {%
0.359800338745117 -0.355319380760193
-0.347306132316589 0.351787686347961
};
\addplot [semithick, white!50.1960784313725!black]
table {%
0.425703644752502 -0.282606363296509
-0.347306132316589 0.351787686347961
};
\addplot [semithick, white!50.1960784313725!black]
table {%
0.425703644752502 -0.282606363296509
-0.405765771865845 0.272964239120483
};
\addplot [semithick, white!50.1960784313725!black]
table {%
0.476155042648315 -0.198433399200439
-0.405765771865845 0.272964239120483
};
\addplot [semithick, white!50.1960784313725!black]
table {%
0.476155042648315 -0.198433399200439
-0.447724223136902 0.184250593185425
};
\addplot [semithick, white!50.1960784313725!black]
table {%
0.476155042648315 -0.198433399200439
-0.480785131454468 0.0918514728546143
};
\addplot [semithick, white!50.1960784313725!black]
table {%
0.5 -0.103239297866821
-0.480785131454468 0.0918514728546143
};
\addplot [semithick, white!50.1960784313725!black]
table {%
0.5 -0.103239297866821
-0.495184659957886 -0.00522196292877197
};
\addplot [semithick, white!50.1960784313725!black]
table {%
0.5 -0.103239297866821
-0.5 -0.103239297866821
};
\addplot [semithick, white!50.1960784313725!black]
table {%
0.495184659957886 -0.00522196292877197
-0.5 -0.103239297866821
};
\addplot [semithick, white!50.1960784313725!black]
table {%
0.480785131454468 0.0918514728546143
-0.5 -0.103239297866821
};
\addplot [very thick, black]
table {%
-0.476155042648315 -0.198433399200439
-0.425703644752502 -0.282606363296509
-0.359800338745117 -0.355319380760193
-0.280977964401245 -0.413777947425842
-0.192265510559082 -0.455735802650452
-0.0970721244812012 -0.479580402374268
-0 -0.493979692459106
0.0970721244812012 -0.479580402374268
0.192265510559082 -0.455735802650452
0.280977964401245 -0.413777947425842
0.359800338745117 -0.355319380760193
0.425703644752502 -0.282606363296509
0.476155042648315 -0.198433399200439
0.5 -0.103239297866821
0.495184659957886 -0.00522196292877197
0.480785131454468 0.0918514728546143
0.447724223136902 0.184250593185425
0.405765771865845 0.272964239120483
0.347306132316589 0.351787686347961
0.274592041969299 0.417691707611084
0.190417766571045 0.468143820762634
0.0980181694030762 0.50120484828949
0 0.506020307540894
-0.0980181694030762 0.50120484828949
-0.190417766571045 0.468143820762634
-0.274592041969299 0.417691707611084
-0.347306132316589 0.351787686347961
-0.405765771865845 0.272964239120483
-0.447724223136902 0.184250593185425
-0.480785131454468 0.0918514728546143
-0.495184659957886 -0.00522196292877197
-0.5 -0.103239297866821
-0.476155042648315 -0.198433399200439
};
\end{axis}

\end{tikzpicture} &
\begin{tikzpicture}

\begin{axis}[
height=\polysize,
hide x axis,
hide y axis,
tick align=outside,
tick pos=left,
width=\polysize,
x grid style={white!69.0196078431373!black},
xmin=-0.549999999980157, xmax=0.550000000019843,
xtick style={color=black},
xtick={-0.6,-0.4,-0.2,0,0.2,0.4,0.6},
xticklabels={\ensuremath{-}0.6,\ensuremath{-}0.4,\ensuremath{-}0.2,0.0,0.2,0.4,0.6},
y grid style={white!69.0196078431373!black},
ymin=-0.565046987010268, ymax=0.534953012989731,
ytick style={color=black},
ytick={-0.6,-0.4,-0.2,0,0.2,0.4,0.6},
yticklabels={\ensuremath{-}0.6,\ensuremath{-}0.4,\ensuremath{-}0.2,0.0,0.2,0.4,0.6}
]
\addplot [semithick, white!50.1960784313725!black]
table {%
-0.283807635307312 0.423778772354126
0.187590479850769 -0.458141684532166
};
\addplot [semithick, white!50.1960784313725!black]
table {%
-0.283807635307312 0.423778772354126
0.27176296710968 -0.407690525054932
};
\addplot [semithick, white!50.1960784313725!black]
table {%
-0.283807635307312 0.423778772354126
0.350584983825684 -0.349232196807861
};
\addplot [semithick, white!50.1960784313725!black]
table {%
-0.356521844863892 0.357874512672424
0.350584983825684 -0.349232196807861
};
\addplot [semithick, white!50.1960784313725!black]
table {%
-0.356521844863892 0.357874512672424
0.41648805141449 -0.276519417762756
};
\addplot [semithick, white!50.1960784313725!black]
table {%
-0.414981603622437 0.279050827026367
0.41648805141449 -0.276519417762756
};
\addplot [semithick, white!50.1960784313725!black]
table {%
-0.414981603622437 0.279050827026367
0.466939330101013 -0.192346572875977
};
\addplot [semithick, white!50.1960784313725!black]
table {%
-0.456940054893494 0.190337061882019
0.466939330101013 -0.192346572875977
};
\addplot [semithick, white!50.1960784313725!black]
table {%
-0.456940054893494 0.190337061882019
0.5 -0.0999484062194824
};
\addplot [semithick, white!50.1960784313725!black]
table {%
-0.480785131454468 0.0951424837112427
0.5 -0.0999484062194824
};
\addplot [semithick, white!50.1960784313725!black]
table {%
-0.495184659957886 -0.00193095207214355
0.5 -0.0999484062194824
};
\addplot [semithick, white!50.1960784313725!black]
table {%
-0.5 -0.0999484062194824
0.5 -0.0999484062194824
};
\addplot [semithick, white!50.1960784313725!black]
table {%
-0.5 -0.0999484062194824
0.495184659957886 -0.00193095207214355
};
\addplot [semithick, white!50.1960784313725!black]
table {%
-0.5 -0.0999484062194824
0.480785131454468 0.0951424837112427
};
\addplot [semithick, white!50.1960784313725!black]
table {%
-0.5 -0.0999484062194824
0.456940054893494 0.190337061882019
};
\addplot [semithick, white!50.1960784313725!black]
table {%
-0.466939330101013 -0.192346572875977
0.456940054893494 0.190337061882019
};
\addplot [semithick, white!50.1960784313725!black]
table {%
-0.466939330101013 -0.192346572875977
0.414981603622437 0.279050827026367
};
\addplot [semithick, white!50.1960784313725!black]
table {%
-0.41648805141449 -0.276519417762756
0.414981603622437 0.279050827026367
};
\addplot [semithick, white!50.1960784313725!black]
table {%
-0.41648805141449 -0.276519417762756
0.356521844863892 0.357874512672424
};
\addplot [semithick, white!50.1960784313725!black]
table {%
-0.350584983825684 -0.349232196807861
0.356521844863892 0.357874512672424
};
\addplot [semithick, white!50.1960784313725!black]
table {%
-0.350584983825684 -0.349232196807861
0.283807635307312 0.423778772354126
};
\addplot [semithick, white!50.1960784313725!black]
table {%
-0.27176296710968 -0.407690525054932
0.283807635307312 0.423778772354126
};
\addplot [semithick, white!50.1960784313725!black]
table {%
-0.187590479850769 -0.458141684532166
0.283807635307312 0.423778772354126
};
\addplot [semithick, white!50.1960784313725!black]
table {%
-0.187590479850769 -0.458141684532166
0.195093154907227 0.465737700462341
};
\addplot [semithick, white!50.1960784313725!black]
table {%
-0.0951930284500122 -0.491202116012573
0.195093154907227 0.465737700462341
};
\addplot [semithick, white!50.1960784313725!black]
table {%
0 -0.515047073364258
0.195093154907227 0.465737700462341
};
\addplot [semithick, white!50.1960784313725!black]
table {%
0 -0.515047073364258
0.0980186462402344 0.480137586593628
};
\addplot [semithick, white!50.1960784313725!black]
table {%
0 -0.515047073364258
-0 0.484953045845032
};
\addplot [semithick, white!50.1960784313725!black]
table {%
0 -0.515047073364258
-0.0980186462402344 0.480137586593628
};
\addplot [semithick, white!50.1960784313725!black]
table {%
0 -0.515047073364258
-0.195093154907227 0.465737700462341
};
\addplot [semithick, white!50.1960784313725!black]
table {%
0.0951930284500122 -0.491202116012573
-0.195093154907227 0.465737700462341
};
\addplot [semithick, white!50.1960784313725!black]
table {%
0.187590479850769 -0.458141684532166
-0.195093154907227 0.465737700462341
};
\addplot [very thick, black]
table {%
-0.283807635307312 0.423778772354126
-0.356521844863892 0.357874512672424
-0.414981603622437 0.279050827026367
-0.456940054893494 0.190337061882019
-0.480785131454468 0.0951424837112427
-0.495184659957886 -0.00193095207214355
-0.5 -0.0999484062194824
-0.466939330101013 -0.192346572875977
-0.41648805141449 -0.276519417762756
-0.350584983825684 -0.349232196807861
-0.27176296710968 -0.407690525054932
-0.187590479850769 -0.458141684532166
-0.0951930284500122 -0.491202116012573
0 -0.515047073364258
0.0951930284500122 -0.491202116012573
0.187590479850769 -0.458141684532166
0.27176296710968 -0.407690525054932
0.350584983825684 -0.349232196807861
0.41648805141449 -0.276519417762756
0.466939330101013 -0.192346572875977
0.5 -0.0999484062194824
0.495184659957886 -0.00193095207214355
0.480785131454468 0.0951424837112427
0.456940054893494 0.190337061882019
0.414981603622437 0.279050827026367
0.356521844863892 0.357874512672424
0.283807635307312 0.423778772354126
0.195093154907227 0.465737700462341
0.0980186462402344 0.480137586593628
-0 0.484953045845032
-0.0980186462402344 0.480137586593628
-0.195093154907227 0.465737700462341
-0.283807635307312 0.423778772354126
};
\end{axis}

\end{tikzpicture} \\
    $\underline{3.140331156954}619\dots$ & 
    $\underline{3.140331156954}543\dots$ & 
    $\underline{3.140331156954}350\dots$   
  \end{tabular}
  \caption{Three small triacontadigons (32-gons) with locally maximal perimeter. The first 13 decimal digits of the perimeter (underlined) coincide.}
  \label{fig:polys_n_32}
\end{figure}

\begin{figure}[htbp]
  \centering
  \setlength\polysize{5.6cm}
\begin{tikzpicture}

\definecolor{darkgray176}{RGB}{176,176,176}
\definecolor{gray}{RGB}{128,128,128}

\begin{axis}[
height=\polysize,
hide x axis,
hide y axis,
tick align=outside,
tick pos=left,
width=\polysize,
x grid style={darkgray176},
xmin=-0.55, xmax=0.55,
xtick style={color=black},
y grid style={darkgray176},
ymin=-0.05, ymax=1.05,
ytick style={color=black}
]
\addplot [semithick, gray]
table {%
0 0
-0 1
};
\addplot [semithick, gray]
table {%
0.19701886177063 0.0196002721786499
-0 1
};
\addplot [semithick, gray]
table {%
0.19701886177063 0.0196002721786499
-0.185775876045227 0.94343376159668
};
\addplot [semithick, gray]
table {%
0.19701886177063 0.0196002721786499
-0.357115745544434 0.852027416229248
};
\addplot [semithick, gray]
table {%
0.34975790977478 0.144687533378601
-0.357115745544434 0.852027416229248
};
\addplot [semithick, gray]
table {%
0.34975790977478 0.144687533378601
-0.480774164199829 0.701658248901367
};
\addplot [semithick, gray]
table {%
0.442771673202515 0.318170309066772
-0.480774164199829 0.701658248901367
};
\addplot [semithick, gray]
table {%
0.5 0.50651216506958
-0.480774164199829 0.701658248901367
};
\addplot [semithick, gray]
table {%
0.5 0.50651216506958
-0.5 0.50651216506958
};
\addplot [semithick, gray]
table {%
0.480774164199829 0.701658248901367
-0.5 0.50651216506958
};
\addplot [semithick, gray]
table {%
0.480774164199829 0.701658248901367
-0.442771673202515 0.318170309066772
};
\addplot [semithick, gray]
table {%
0.480774164199829 0.701658248901367
-0.34975790977478 0.144687533378601
};
\addplot [semithick, gray]
table {%
0.357115745544434 0.852027416229248
-0.34975790977478 0.144687533378601
};
\addplot [semithick, gray]
table {%
0.357115745544434 0.852027416229248
-0.19701886177063 0.0196002721786499
};
\addplot [semithick, gray]
table {%
0.185775876045227 0.94343376159668
-0.19701886177063 0.0196002721786499
};
\addplot [semithick, gray]
table {%
-0 1
-0.19701886177063 0.0196002721786499
};
\addplot [very thick, black]
table {%
0 0
0.19701886177063 0.0196002721786499
0.34975790977478 0.144687533378601
0.442771673202515 0.318170309066772
0.5 0.50651216506958
0.480774164199829 0.701658248901367
0.357115745544434 0.852027416229248
0.185775876045227 0.94343376159668
-0 1
-0.185775876045227 0.94343376159668
-0.357115745544434 0.852027416229248
-0.480774164199829 0.701658248901367
-0.5 0.50651216506958
-0.442771673202515 0.318170309066772
-0.34975790977478 0.144687533378601
-0.19701886177063 0.0196002721786499
0 0
};
\end{axis}

\end{tikzpicture}
\begin{tikzpicture}

\definecolor{darkgray176}{RGB}{176,176,176}
\definecolor{gray}{RGB}{128,128,128}

\begin{axis}[
height=\polysize,
hide x axis,
hide y axis,
tick align=outside,
tick pos=left,
width=\polysize,
x grid style={darkgray176},
xmin=-0.55, xmax=0.55,
xtick style={color=black},
y grid style={darkgray176},
ymin=-0.05, ymax=1.05,
ytick style={color=black}
]
\addplot [semithick, gray]
table {%
0 0
0.0490677356719971 0.998795509338379
};
\addplot [semithick, gray]
table {%
0 0
0 1
};
\addplot [semithick, gray]
table {%
0 0
-0.0490677356719971 0.998795509338379
};
\addplot [semithick, gray]
table {%
0.0489494800567627 0.00361073017120361
-0.0490677356719971 0.998795509338379
};
\addplot [semithick, gray]
table {%
0.0976628065109253 0.00961899757385254
-0.0490677356719971 0.998795509338379
};
\addplot [semithick, gray]
table {%
0.14602267742157 0.018010139465332
-0.0490677356719971 0.998795509338379
};
\addplot [semithick, gray]
table {%
0.193912506103516 0.0287642478942871
-0.0490677356719971 0.998795509338379
};
\addplot [semithick, gray]
table {%
0.241217017173767 0.0418550968170166
-0.0490677356719971 0.998795509338379
};
\addplot [semithick, gray]
table {%
0.287822246551514 0.0572513341903687
-0.0490677356719971 0.998795509338379
};
\addplot [semithick, gray]
table {%
0.287822246551514 0.0572513341903687
-0.0948612689971924 0.981130838394165
};
\addplot [semithick, gray]
table {%
0.287822246551514 0.0572513341903687
-0.139732956886292 0.961240768432617
};
\addplot [semithick, gray]
table {%
0.287822246551514 0.0572513341903687
-0.183574557304382 0.939172744750977
};
\addplot [semithick, gray]
table {%
0.287822246551514 0.0572513341903687
-0.226280570030212 0.914979934692383
};
\addplot [semithick, gray]
table {%
0.287822246551514 0.0572513341903687
-0.267748117446899 0.888720989227295
};
\addplot [semithick, gray]
table {%
0.327951192855835 0.0855134725570679
-0.267748117446899 0.888720989227295
};
\addplot [semithick, gray]
table {%
0.327951192855835 0.0855134725570679
-0.306442022323608 0.858523845672607
};
\addplot [semithick, gray]
table {%
0.327951192855835 0.0855134725570679
-0.343607664108276 0.826464653015137
};
\addplot [semithick, gray]
table {%
0.363499045372009 0.119357824325562
-0.343607664108276 0.826464653015137
};
\addplot [semithick, gray]
table {%
0.363499045372009 0.119357824325562
-0.37745201587677 0.790916800498962
};
\addplot [semithick, gray]
table {%
0.39555835723877 0.156523466110229
-0.37745201587677 0.790916800498962
};
\addplot [semithick, gray]
table {%
0.425755500793457 0.195217490196228
-0.37745201587677 0.790916800498962
};
\addplot [semithick, gray]
table {%
0.454017519950867 0.235346555709839
-0.37745201587677 0.790916800498962
};
\addplot [semithick, gray]
table {%
0.454017519950867 0.235346555709839
-0.403711080551147 0.749449253082275
};
\addplot [semithick, gray]
table {%
0.454017519950867 0.235346555709839
-0.427903652191162 0.706743240356445
};
\addplot [semithick, gray]
table {%
0.476085662841797 0.27918815612793
-0.427903652191162 0.706743240356445
};
\addplot [semithick, gray]
table {%
0.476085662841797 0.27918815612793
-0.447793960571289 0.661871671676636
};
\addplot [semithick, gray]
table {%
0.476085662841797 0.27918815612793
-0.465458512306213 0.616078019142151
};
\addplot [semithick, gray]
table {%
0.476085662841797 0.27918815612793
-0.480854749679565 0.569472789764404
};
\addplot [semithick, gray]
table {%
0.489176511764526 0.326492667198181
-0.480854749679565 0.569472789764404
};
\addplot [semithick, gray]
table {%
0.489176511764526 0.326492667198181
-0.491608738899231 0.521583080291748
};
\addplot [semithick, gray]
table {%
0.489176511764526 0.326492667198181
-0.5 0.473223209381104
};
\addplot [semithick, gray]
table {%
0.495184659957886 0.375205993652344
-0.5 0.473223209381104
};
\addplot [semithick, gray]
table {%
0.498795509338379 0.424155473709106
-0.5 0.473223209381104
};
\addplot [semithick, gray]
table {%
0.5 0.473223209381104
-0.5 0.473223209381104
};
\addplot [semithick, gray]
table {%
0.5 0.473223209381104
-0.498795509338379 0.424155473709106
};
\addplot [semithick, gray]
table {%
0.5 0.473223209381104
-0.495184659957886 0.375205993652344
};
\addplot [semithick, gray]
table {%
0.5 0.473223209381104
-0.489176511764526 0.326492667198181
};
\addplot [semithick, gray]
table {%
0.491608738899231 0.521583080291748
-0.489176511764526 0.326492667198181
};
\addplot [semithick, gray]
table {%
0.480854749679565 0.569472789764404
-0.489176511764526 0.326492667198181
};
\addplot [semithick, gray]
table {%
0.480854749679565 0.569472789764404
-0.476085662841797 0.27918815612793
};
\addplot [semithick, gray]
table {%
0.465458512306213 0.616078019142151
-0.476085662841797 0.27918815612793
};
\addplot [semithick, gray]
table {%
0.447793960571289 0.661871671676636
-0.476085662841797 0.27918815612793
};
\addplot [semithick, gray]
table {%
0.427903652191162 0.706743240356445
-0.476085662841797 0.27918815612793
};
\addplot [semithick, gray]
table {%
0.427903652191162 0.706743240356445
-0.454017519950867 0.235346555709839
};
\addplot [semithick, gray]
table {%
0.403711080551147 0.749449253082275
-0.454017519950867 0.235346555709839
};
\addplot [semithick, gray]
table {%
0.37745201587677 0.790916800498962
-0.454017519950867 0.235346555709839
};
\addplot [semithick, gray]
table {%
0.37745201587677 0.790916800498962
-0.425755500793457 0.195217490196228
};
\addplot [semithick, gray]
table {%
0.37745201587677 0.790916800498962
-0.39555835723877 0.156523466110229
};
\addplot [semithick, gray]
table {%
0.37745201587677 0.790916800498962
-0.363499045372009 0.119357824325562
};
\addplot [semithick, gray]
table {%
0.343607664108276 0.826464653015137
-0.363499045372009 0.119357824325562
};
\addplot [semithick, gray]
table {%
0.343607664108276 0.826464653015137
-0.327951192855835 0.0855134725570679
};
\addplot [semithick, gray]
table {%
0.306442022323608 0.858523845672607
-0.327951192855835 0.0855134725570679
};
\addplot [semithick, gray]
table {%
0.267748117446899 0.888720989227295
-0.327951192855835 0.0855134725570679
};
\addplot [semithick, gray]
table {%
0.267748117446899 0.888720989227295
-0.287822246551514 0.0572513341903687
};
\addplot [semithick, gray]
table {%
0.226280570030212 0.914979934692383
-0.287822246551514 0.0572513341903687
};
\addplot [semithick, gray]
table {%
0.183574557304382 0.939172744750977
-0.287822246551514 0.0572513341903687
};
\addplot [semithick, gray]
table {%
0.139732956886292 0.961240768432617
-0.287822246551514 0.0572513341903687
};
\addplot [semithick, gray]
table {%
0.0948612689971924 0.981130838394165
-0.287822246551514 0.0572513341903687
};
\addplot [semithick, gray]
table {%
0.0490677356719971 0.998795509338379
-0.287822246551514 0.0572513341903687
};
\addplot [semithick, gray]
table {%
0.0490677356719971 0.998795509338379
-0.241217017173767 0.0418550968170166
};
\addplot [semithick, gray]
table {%
0.0490677356719971 0.998795509338379
-0.193912506103516 0.0287642478942871
};
\addplot [semithick, gray]
table {%
0.0490677356719971 0.998795509338379
-0.14602267742157 0.018010139465332
};
\addplot [semithick, gray]
table {%
0.0490677356719971 0.998795509338379
-0.0976628065109253 0.00961899757385254
};
\addplot [semithick, gray]
table {%
0.0490677356719971 0.998795509338379
-0.0489494800567627 0.00361073017120361
};
\addplot [very thick, black]
table {%
0 0
0.0489494800567627 0.00361073017120361
0.0976628065109253 0.00961899757385254
0.14602267742157 0.018010139465332
0.193912506103516 0.0287642478942871
0.241217017173767 0.0418550968170166
0.287822246551514 0.0572513341903687
0.327951192855835 0.0855134725570679
0.363499045372009 0.119357824325562
0.39555835723877 0.156523466110229
0.425755500793457 0.195217490196228
0.454017519950867 0.235346555709839
0.476085662841797 0.27918815612793
0.489176511764526 0.326492667198181
0.495184659957886 0.375205993652344
0.498795509338379 0.424155473709106
0.5 0.473223209381104
0.491608738899231 0.521583080291748
0.480854749679565 0.569472789764404
0.465458512306213 0.616078019142151
0.447793960571289 0.661871671676636
0.427903652191162 0.706743240356445
0.403711080551147 0.749449253082275
0.37745201587677 0.790916800498962
0.343607664108276 0.826464653015137
0.306442022323608 0.858523845672607
0.267748117446899 0.888720989227295
0.226280570030212 0.914979934692383
0.183574557304382 0.939172744750977
0.139732956886292 0.961240768432617
0.0948612689971924 0.981130838394165
0.0490677356719971 0.998795509338379
0 1
-0.0490677356719971 0.998795509338379
-0.0948612689971924 0.981130838394165
-0.139732956886292 0.961240768432617
-0.183574557304382 0.939172744750977
-0.226280570030212 0.914979934692383
-0.267748117446899 0.888720989227295
-0.306442022323608 0.858523845672607
-0.343607664108276 0.826464653015137
-0.37745201587677 0.790916800498962
-0.403711080551147 0.749449253082275
-0.427903652191162 0.706743240356445
-0.447793960571289 0.661871671676636
-0.465458512306213 0.616078019142151
-0.480854749679565 0.569472789764404
-0.491608738899231 0.521583080291748
-0.5 0.473223209381104
-0.498795509338379 0.424155473709106
-0.495184659957886 0.375205993652344
-0.489176511764526 0.326492667198181
-0.476085662841797 0.27918815612793
-0.454017519950867 0.235346555709839
-0.425755500793457 0.195217490196228
-0.39555835723877 0.156523466110229
-0.363499045372009 0.119357824325562
-0.327951192855835 0.0855134725570679
-0.287822246551514 0.0572513341903687
-0.241217017173767 0.0418550968170166
-0.193912506103516 0.0287642478942871
-0.14602267742157 0.018010139465332
-0.0976628065109253 0.00961899757385254
-0.0489494800567627 0.00361073017120361
0 0
};
\end{axis}

\end{tikzpicture}
  \input{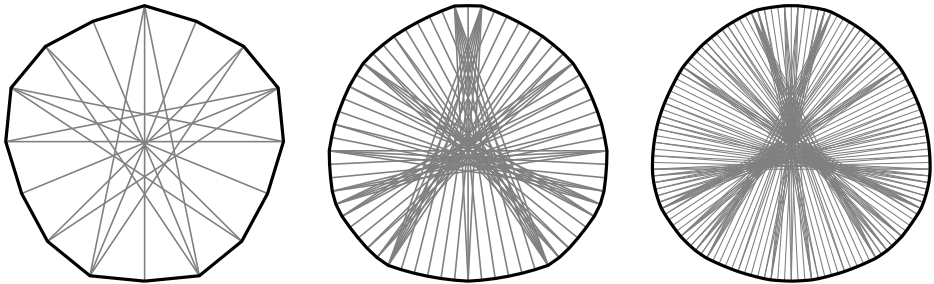}
  \caption{Current perimeter records of small $n$-gons for $n=2^s$, $s=4, 6, 7$ (with $s=5$ being depicted in Fig.~\ref{fig:polys_n_32}, left).}
  \label{fig:polys_vary_n}
\end{figure}

\begin{table}[htbp]
  \centering
  \begin{tabular}{cccc}
    $n$ & $\pi^9 / (8 n^8)$ & $\bar{u}_n - p(P_n)$ & Perimeter $p(P_n)$ \\
    \hline
    4   & $5.686 \cdot 10^{-2}$ & $2.619 \cdot 10^{-2}$   & 3.03527618041008304939559535049619331339627560527972205 \\ &&& $\hookrightarrow$ 52560128292602278989952079876894718987769986620\dots \\ 
    8   & $2.221 \cdot 10^{-4}$ & $2.980 \cdot 10^{-4}$   & 3.12114713405983135386465950363808653090954216646976012 \\ &&& $\hookrightarrow$ 24524789123816403490428894959252350355455226792\dots \\ 
    16  & $8.676 \cdot 10^{-7}$ & $7.741 \cdot 10^{-7}$   & 3.13654771648660738608596703194122822729813676580923269 \\ &&& $\hookrightarrow$ 27892182035777457554738176289058573625428211593\dots \\ 
    32  & $3.389 \cdot 10^{-9}$ & $1.335 \cdot 10^{-13}$  & 3.14033115695461936582540138057745867231205309833952186 \\ &&& $\hookrightarrow$ 99104148559468837774634543964164383698560055119\dots \\ 
    64  & $1.324 \cdot 10^{-11}$ & $2.836 \cdot 10^{-23}$ & 3.14127725093277286806199141550246829795626209630809641 \\ &&& $\hookrightarrow$ 11750773439718362183509788657317267672710085186\dots \\ 
    128 & $5.171 \cdot 10^{-14}$ & $1.816 \cdot 10^{-38}$ & 3.14151380114430107632851505945682230791714977539831260 \\ &&& $\hookrightarrow$ 12200604676901080305902623648703203853047686174\dots \\ 
    \hline
  \end{tabular}
  \caption{Perimeters of the polygons of Fig.~\ref{fig:polys_vary_n}, their distances to Reinhardt's upper bound, and the asymptotic $\pi^9 / (8 n^8)$. The known solution $2+4\sin\frac{\pi}{12}$ for $n=4$ is provided for comparison.}
  \label{tab:perimeters_vary_n}
\end{table}

\subsection{Structure}

In Sec.~\eqref{sec:zonogon}, we use a zonogon construction to derive a MINLP for the perimeter problem. We discuss the combinatorics of the problem in Sec.~\ref{sec:code_combinatorics} and propose a two-phase algorithm to compute polygons with large perimeter to arbitrary accuracy in Sec.~\ref{sec:two_phase}. We present our computational experience with the two-phase approach in Sec.~\ref{sec:numerics}. Finally, we show that the maximum perimeter for the octagon is an algebraic number in Sec.~\ref{sec:perim_algebraic} and provide the corresponding integer polygon in the Appendix.

\section{MINLP derivation using zonogons} \label{sec:zonogon}

\begin{figure}[btp]
  \centering
  \begin{tikzpicture}[scale=3,line cap=round,line join=round,x=1cm,y=1cm]
    \draw[very thick,dashed] (-0.44442725609268774,0.8958149440827188) -- (0.5912742926107506,0.8064705269863596) -- (0.6452255772671045,0.4164170631939196) -- (0,0) -- cycle;
    \draw[very thick,dashed] (0,0) -- (1.0357015487034382,-0.08934441709635921) -- (1.0896528333597923,-0.4793978808887992) -- (0.44442725609268774,-0.8958149440827188) -- cycle;
    \draw[very thick,dashed] (-1.0357015487034382,0.08934441709635921) -- (0,0) -- (0.053951284656353904,-0.39005346379244) -- (-0.5912742926107506,-0.8064705269863596) -- cycle;
    \draw[very thick,dashed] (-1.0896528333597923,0.4793978808887992) -- (-0.053951284656353904,0.39005346379244) -- (0,0) -- (-0.6452255772671045,-0.4164170631939196) -- cycle;
    \draw (0.6,-0.35) node {$P_1$};
    \draw (0.2,0.5) node {$P_2$};
    \draw (-0.4,0.2) node {$P_3$};
    \draw (-0.2,-0.4) node {$P_4$};
    \draw 
      (-0.44442725609268774,0.8958149440827188) node[above] {$z_1 = v_1 - v_2$} -- 
      (-1.0896528333597923,0.4793978808887992) node[left] {$z_2 = v_1 - v_3$} --
      (-1.0357015487034382,0.08934441709635921) node[left] {$z_3 = v_1 - v_4$} --
      (-0.5912742926107506,-0.8064705269863596) node[left] {$z_4 = v_2 - v_4$} --
      (0.44442725609268774,-0.8958149440827188) node[below] {$z_5 = v_2 - v_1$} --
      (1.0896528333597923,-0.4793978808887992) node[right] {$z_6 = v_3 - v_1$} --
      (1.0357015487034382,-0.08934441709635921) node[right] {$z_7 = v_4 - v_1$} --
      (0.5912742926107506,0.8064705269863596) node[right] {$z_8 = v_4 - v_2$} -- cycle;
    \end{tikzpicture}
  \caption{Construction of the zonogon $Z$ (solid edges) from the shifts 
  $P_1$, $P_2$, $P_3$, $P_4$ (dashed edges) of the polygon $P$.}
  \label{fig:zonogon}
\end{figure}
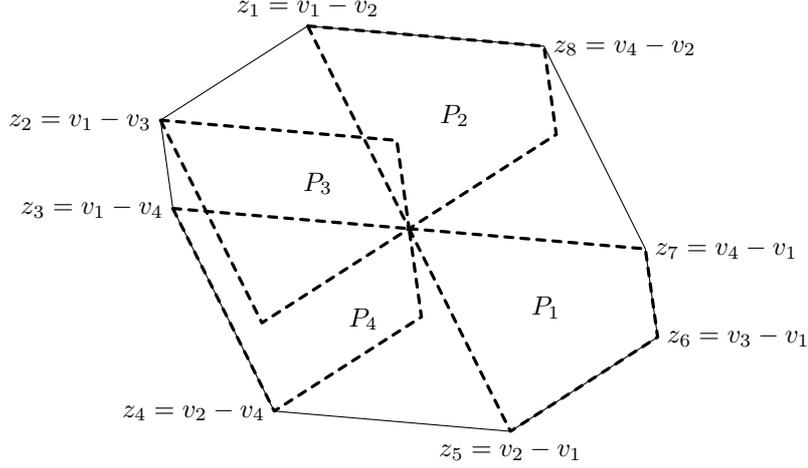

Let $P$ be a planar convex polygon with $n \ge 3$ vertices $v_1, \dotsc, v_{n} \in \mathbb{C}$, enumerated in counterclockwise order. We interpret $P$ here as a closed subset of $\mathbb{C}$ in the sense of the convex hull of its vertices $P = \cvx \{ v_1, \dotsc, v_{n}\}$. Let us consider $n$ shifts of $P$ defined by
\[
  P_{l} = P - v_{l}, \quad l=1,\dotsc,n,
\]
i.e., the $l$-th vertex of $P$ is moved to the origin. The convex hull 
(compare Fig.~\ref{fig:zonogon})
\begin{equation} \label{eqn:zonogon_construction}
  Z 
  = \left\{ \sum_{k,l=1}^{n} \lambda_{kl}(v_{k} - v_{l}) \mid \sum_{k,l=1}^{n} \lambda_{kl} = 1, \; \lambda_{kl} \ge 0, \; k, l = 1, \dotsc, n \right\}
\end{equation}
of the shifted polygons $P_1, \dotsc, P_{n}$ is the \emph{zonogon} generated by $P$.
Indeed, if $z \in Z$, then $-z \in Z$ by transposition of $\lambda$, which means that $Z$ is
centrally symmetric with respect to the origin.

\begin{lemma} \label{lem:B}
  The zonogon $Z$ is contained in the unit disc if and only if $P$ is small. 
\end{lemma}
\begin{proof}
  Let $\abs{z} \le 1$ for all $z \in Z$. Choosing $\lambda_{kl} = 1$ (and all other components zero) delivers 
  $v_k-v_l \in Z$, hence $\abs{v_{k} - v_{l}}\le1$ for all $k, l = 1, \dotsc, n$. Hence, $P$ is small.

  Now assume $\abs{v_{k} - v_{l}} \le 1$ for all $k, l = 1, \dotsc, n$. Then, $\abs{z} = \abs{\smash[b]{\sum_{k,l=1}^{n} \lambda_{kl}(v_{k} - v_{l})}} \le \sum_{k,l=1}^{n} \lambda_{kl} = 1.$
\end{proof}

Because $P$ is a convex polygon, the zonogon $Z$ is polygonal as well. 
By the representation \eqref{eqn:zonogon_construction} of $Z$ as convex hull of the points
$v_k-v_l$, its vertices have to be of this form.
Choosing a vertex of $Z$ as a start, we traverse the vertices counterclockwise to get
\[
  z_{j} = v_{k_{j}} - v_{l_{j}},\quad k_j,l_j \in\{1,\dotsc,n\}, \quad j=1,\dotsc,2n.
\]
The index tuples $(k_j)$ and $(l_j)$ are such that for each $j$ either
\begin{equation} \label{eqn:code_one}
  k_{j+1} = k_{j} + 1 \quad \text{and} \quad l_{j+1} = l_{j},
\end{equation}
which means a move along the side $v_{k_j}v_{{k_j}+1}$ of $P_{l_{j}}$,
or
\begin{equation} \label{eqn:code_minus_one}
  k_{j+1} = k_{j} \quad \text{and} \quad l_{j+1} = l_{j} + 1,
\end{equation}
which means a move to the vertex $v_{k_j}$ of the next shift $P_{l_{j} + 1}$.%
\footnote{All indices $k_{j}$, $l_{j}$ are understood as members of the residue system $\{1,\dotsc,n\}$
modulo $n$, which we shall silently assume for all upcoming index computations.}
After visiting all vertices of $Z$ we return to $z_1$. Then the indices $k_{j}$ and $l_{j}$ 
have been incremented $n$ times each, which yields that $Z$ has $2n$ vertices. We will assume that
$P$ does not have parallel sides, which implies that the vertices of $Z$ are nondegenerate.

\begin{definition} \label{def:code}
  Let $P$ be a convex polygon such that the generated zonogon $Z$ has $2n$ vertices. We define the \emph{code} of $P$ as $c=(c_1,\dotsc,c_{2n}) \in \{\pm 1\}^{2n}$, where $c_j := 1$ if~\eqref{eqn:code_one} holds
respectively $c_j := -1$ if~\eqref{eqn:code_minus_one} holds.
\end{definition}
The dependence of the code on the choice of the start vertex is not problematic, because we always consider equivalence classes of codes, where two codes $c$ and $c'$ are equivalent if there is an element of the dihedral group $\sigma \in D_{2n}$ such that $c_{\sigma(j)} = c'_j$ for $j = 1, \dotsc, 2n$,
see~\autoref{sec:code_combinatorics}.
It is most convenient to assign the signs $c_j$ to the sides $z_{j+1}-z_j$ of $Z$.
We remark that $c_{n+j} = -c_j$ for $j = 1, \dotsc, n$, because exactly one of the two opposite sides of $Z$ is a side of one of the shifts of $P$. The code corresponding to Fig.~\ref{fig:zonogon} is $c = (--+-++-\,+)$, where the signs in the shortened notation have to be expanded to $\pm 1$. The corresponding index tuples are
\[
  k = (1,1,1,2,2,3,4,4) \quad \text{and} \quad
  l = (2,3,4,4,1,1,1,2).
\]
Note that $k$ and $l$, considered as circular tuples, are just rotations of each other.

We now parametrize the vertices of $Z$ in polar form as
\[
  z_j = d_j \exp(\ii \varphi_j) \quad \text{for } j = 1, \dotsc, 2n.
\]
The central symmetry of $Z$ dictates that $z_{n+j}=-z_{j}$, i.e. $d_{n+j}=d_j$ and
$\varphi_{n+j} = \varphi_j + \pi$ (modulo $2 \pi$). 
The side vectors of $P$ are written as
\[
	s_j = \begin{cases} v_{k_{j+1}}-v_{k_j}, & c_j=1, \\   v_{l_{j}}-v_{l_{j+1}}, &c_j=-1, \end{cases} \Bigg\}
	=c_j(z_{j+1}-z_{j}) =: c_j r_j \exp(\ii \rho_j) \quad \text{for } j = 1, \dotsc, n.
\]
Hence, the construction of the generated zonogon $Z$ implies that $s_{j} = c_j(z_{j+1}-z_j)$ for $j=1,\dotsc,n$.
We explicitly need to enforce that
\begin{equation} \label{eqn:complex_closed_edges}
  0 = \sum_{j=1}^{n} s_j = \sum_{j=1}^{n} c_j (z_{j+1} - z_{j}),
\end{equation}
so that the boundary curve of $P$ is closed. 

We are now ready to formulate the perimeter problem as the MINLP
\begin{equation} \label{eqn:minlp_general}
  \begin{alignedat}{3}
    \max~ & \sum_{j=1}^n r_j \\
    \text{s.t.~}
    & 0\le r_j, \quad 0\le d_j\le1, \\
    &r_{j} \exp(\ii \rho_{j})= d_{j+1}\exp(\ii \varphi_{j+1})-d_j \exp(\ii \varphi_j),\quad j=1,\dotsc,n, \\
    & \sum_{j=1}^{n} c_j  r_{j} \exp(\ii \rho_{j})= 0, \\
    & \varphi_1 \le \varphi_2 \le \dotsb \le \varphi_{n+1} = \varphi_1 + \pi, \quad
    \rho_1 \le \rho_2 \le \dotsb \le \rho_{n+1} = \rho_1 + \pi,\\
    & c_j \in \{ \pm 1 \}, \quad j = 1, \dotsc, n.
  \end{alignedat}
\end{equation}

We solved all NLPs that result from MINLP~\eqref{eqn:minlp_general} for each possible code $c$ for $n=4, 8, 16, 32$ numerically and found that $d_j = 1$ always held.
Based on this numerical evidence, we now consider the case that $P$ has a full set of $n$ distinct diameters $v_kv_l$ with $\abs{v_{k} - v_{l}} = 1$ 
(more are not possible, see~\cite{erdos1946sets}). This happens exactly if the zonogon $Z$
is inscribed into the unit circle, i.e., all its vertices $z_{j}$, $j=1, \dotsc, 2n$, lie on the unit circle. 
%
Then the problem allows for the following formulation as a MINLP
\begin{equation} \label{eqn:minlp}
  \begin{alignedat}{3}
    \max~ & \sum_{j=1}^n 2 \sin \frac{\varphi_{j+1} - \varphi_{j}}{2} \\
    \text{s.t.~}
    & \sum_{j=1}^{n} c_j (\cos \varphi_{j+1} - \cos \varphi_{j}) = 0, \quad
    \sum_{j=1}^{n} c_j (\sin \varphi_{j+1} - \sin \varphi_{j}) = 0, \\
    &\varphi_1 \le \varphi_2 \le \dotsb \le \varphi_{n+1} = \varphi_1 + \pi, \\
    & c_j \in \{ \pm 1 \}, \quad j = 1, \dotsc, n,
  \end{alignedat}
\end{equation}
which comprises $n$ continuous real variables $\varphi_1, \dotsc, \varphi_n$ and $n$ ``binary'' variables $c_1, \dotsc, c_n$. The two equality constraints, which are the real and imaginary parts of~\eqref{eqn:complex_closed_edges}, are nonlinear in $\varphi$.
The continuous rotational symmetry of MINLP~\eqref{eqn:minlp} can be eliminated by fixing one of the angles, e.g., $\varphi_1 = 0$.

The MINLP~\eqref{eqn:minlp} is more general than the NLP proposed in~\cite{bingane_maximal_2023}, which is formulated for fixed compositions instead of codes that are still free for optimization. Furthermore, working with absolute angles instead of angle differences yields computationally more favorable sparsity structures, which we shall exploit in Sec.~\ref{sec:two_phase}.

\section{Code combinatorics} \label{sec:code_combinatorics}

The codes $c$ characterize the combinatorial structure of $P$ and $Z$. Clearly, the perimeter problem
allows for congruence transformations. So we only have to deal with the equivalence classes of codes
with respect to the automorphism group of the cycle graph $C_n$, which is known to be isomorphic
to the dihedral group $D_n$.
To enumerate the essentially distinct codes for the MINLPs~\eqref{eqn:minlp_general} and~\eqref{eqn:minlp} we propose two approaches, which are based on the zonogon and the polygon viewpoints. 

\subsection{Zonogon view} 

We consider again the zonogon $Z$ with the code $c \in \{ \pm 1 \}^{2n}$ attached to its sides. Because $c_{n+j} = -c_{j}$, the equivalent codes can be considered as a self-dual two-colored (plus/minus) bracelet (turnover necklace) with $2n$ beads~\cite{palmer1984enumeration} with at most $n-2$ consecutive beads of the same color. This restriction is caused by the fact that at most $n-2$ sides of each $P_l$ may be also sides of $Z$, because $v_l$
has been moved to the origin, which is inside $Z$.

There is exactly one self-dual bracelet of length $2n$ with $n$ consecutive beads: 
\[ (+ + \dotsb + - - \dotsb -).\]
Moreover, self-dual bracelets of length $2n$ with $n-1$ consecutive beads do not exist, because they would need to have the form
\[
  (\underbrace{+ \dotsb +}_{n-1 \text{ times}} - \underbrace{- \dotsb -}_{n-1 \text{ times}} +),
\]
which has $n$ consecutive beads.

Hence, the number of essentially distinct codes is one less than the number of self-dual two-colored bracelets with $2n$ beads, which is known as sequence A007147 in the OEIS~\cite{OEISA007147}.

\subsection{Vertex-oriented polygon view}
Alternatively, we can enumerate codes by looking at the (undirected) diameter graph of the polygon $P$, which consists of the nodes $\{ v_1, \dotsc, v_n \}$ and has edges connecting $v_k$ and $v_l$ if $\abs{v_k - v_l} = 1$. 
Supposing a full set of $n$ diameters, the edges of the diameter graph are in one-to-one correspondence with the code $c$ through the sequences $(l_j)$ and $(k_j)$ because $z_j = v_{k_j} - v_{l_j}$.
Coloring in red each vertex that is incident to only one edge (i.e., appears only once in each of the sequences $(l_j)$ and $(k_j)$) and all other vertices white, we obtain a graph with two-colored nodes. After removal of all red nodes, the remaining part of the diameter graph must be an odd cycle (compare~\cite{erdos1946sets}). The zonogon construction informs us that there can be at most $n-3$ red nodes and, hence, there must be at least three white nodes.
In Fig.~\ref{fig:zonogon}, for instance, the only red vertex would be $v_3$ because $l_6 = 3$ and $l_j \neq 3$ for all $j \neq 6$.
Consequently, the number of distinct codes is equal to the number of bracelets with $n$ beads colored white or red with an odd number of at least three white beads. This sequence is known as sequence A263768 in the OEIS~\cite{OEISA263768}.

As a consequence, we have just proved that
\[
  \text{A263768}(n) = \text{A007147}(n) - 1.
\]
In Tab.~\ref{tab:number_of_codes} we see that although the number of distinct codes grows rapidly, we can completely enumerate all possible codes in a reasonable amount of time at least for $n \le 32$ using fast algorithms~\cite{sawada2001generating}.


\subsection{Integer composition view}
To encode diameter graphs, it is customary to use integer compositions of $n$, which are $m$-tuples of positive integers that sum up to $n$. It is not difficult to derive that there is a one-to-one correspondence to the codes $c \in \{\pm 1\}^{2n}$, by counting the number of consecutive negative entries or, equivalently, the number of consecutive positive entries.

As an example, the code $c = (+--+-++--++-+--\,+)$ of the optimal octagon leads to the composition $(2, 1, 2, 1, 2)$ by counting the consecutive minus signs. When counting consecutive signs, one has to consider the code
as cyclic. Counting plus signs in the example leads to the composition $(1, 2, 2, 1, 2)$ or $(2, 1, 2, 2, 1)$, depending on whether the first plus is wrapped to the end or vice versa. All three compositions are equivalent.
To recover a code from, e.g., the composition $(2, 1, 2, 2, 1)$, we first write down the consecutive groups of negative signs and then fill the blanks, starting left of the first minus group with plus signs going through the composition from its center entry $(2, 2, 1, 2, 1)$
\[
  (-- \quad - \quad -- \quad -- \quad -)
  \quad \leadsto \quad
  (++ -- ++ - + -- ++ -- +\,-),
\]
which is equivalent to the original $c$ under the action of the dihedral group $D_{16}$, in this particular case a shift by five to the right.


\section{Two-phase algorithm using arbitrary precision arithmetic} \label{sec:two_phase}

We propose a two-phase algorithm that computes a code first (Phase~I) and then computes a (local) solution of the fixed code NLP with an arbitrary precision Lagrange-Newton method (Phase~II).

The approach is motivated by taking the regular $2n$-gon for $Z$, i.e., choosing
\[
  \varphi^\ast_j = \frac{j \pi}{n} \quad \text{for } j = 1, \dotsc, 2n.
\]
Inserting $\varphi^\ast$ into MINLP~\eqref{eqn:minlp}, we immediately obtain that the objective coincides with Reinhardt's upper bound $\bar{u}_n$. The idea of Phase~I is now to find a code $c^\ast \in \{\pm 1\}$ with minimal violation of the equality constraints, while in Phase~II we fix $c = c^\ast$ and use variations of a local Lagrange-Newton method with initial guess $(c^\ast, \varphi^\ast)$ (disregarding the angle inequalities).

\subsection{Phase~I: Subset-Sum Problem}
Using $\varphi^\ast$, we obtain for $j = 1, \dotsc, 2n$ that
\[
  z^\ast_j = \exp(\ii \varphi^\ast_j) = \xi^j \quad \text{with } \xi = \exp(\ii \pi / n).
\]
The complex-valued version~\eqref{eqn:complex_closed_edges} of the equality constraint function in MINLP~\eqref{eqn:minlp} then becomes
\begin{equation*} 
  r(c) := \sum_{j=1}^{n} c_j (z^\ast_{j+1} - z^\ast_j)
  = (\xi - 1) \sum_{j=1}^{n} c_j \xi^j.
\end{equation*}
Hence, minimizing the constraint violation means to pick $n$ non-opposing $2n$-th roots of unity such that their sum is as small as possible.

\subsubsection{Case of odd divisors}
We immediately observe that if $n$ has an odd divisor $d \ge 3$, there are (possibly multiple) codes $c^\ast$ for which $r(c^\ast)$ vanishes. In this case $(c^\ast, \varphi^\ast)$ is already a global solution of MINLP~\eqref{eqn:minlp} because it attains the upper bound $\bar{u}_n$. One simple construction for such a code is
\[
  c^\ast_j = (-1)^{\ceil{\frac{j d}{n}}} \quad \text{for } j = 1, \dotsc, 2n,
\]
which satisfies
\[
  c^\ast_{n+j} = (-1)^{d + \ceil{\frac{j d}{n}}} = -c^\ast_{j}
\]
and with the representation $j = p + q \frac{n}{d}$ for $p=1, \dotsc, \frac{n}{d}$, $q = 0, \dotsc, d-1$, also
\begin{align*}
  \frac{r(c^\ast)}{\xi-1} &= \sum_{j=1}^{n} c^\ast_j \xi^j
  = \sum_{q=0}^{d-1} \sum_{p=1}^{n/d} (-1)^{q + \ceil{\frac{pd}{n}}} \xi^{p + \frac{qn}{d}}
  = \left(\sum_{q=0}^{d-1} \left(-\xi^{\frac{n}{d}}\right)^q \right) \left(-\sum_{p=1}^{n/d} \xi^p \right) = 0,
\end{align*}
because the first factor is the sum of all $d$-th roots of unity.
If $n$ is odd, an even simpler choice is $c^\ast_j = (-1)^j$.

\subsubsection{Case of $n$ a power of two} \label{sec:phase_one_power_of_two}
For $n = 2^s$, there is in general no code $c^\ast$ with $r(c^\ast) = 0$ and the problem becomes an interesting variation on a challenging question about sums of roots of unity~\cite{myerson1986unsolved,christie_classifying_2020}.

Exploiting the dihedral symmetry of $c^\ast$ inherent in MINLP~\eqref{eqn:minlp}, we can assume without loss of generality that $\arg r(c^\ast) \in [0, \frac{1}{2}\varphi^\ast_1]$ and minimize the real part of $r(c^\ast)$ (or do the same for $r(c^\ast)/(\xi-1)$).

The computations become easier if we assume Mossinghoff's conjecture to be true: We can then resort to symmetric codes
\begin{equation} \label{eqn:sym_code}
  c_{n-j+1} = -c_j \quad \text{for } j = 1, \dotsc, \tfrac{n}{2}.
\end{equation}
Besides halving the degrees of freedom, this also leads to
\begin{align*}
  \xi^{-1} r(c) &= (1 - \xi^{-1}) \sum_{j=1}^{n} c_j \xi^j
  = \left( \xi^{\frac{1}{2}} - \xi^{-\frac{1}{2}}\right) \sum_{j=1}^{n} c_j \xi^{j-\frac{1}{2}}
  = 4 \ii \Im\left(\xi^{\frac{1}{2}}\right) \sum_{j=1}^{n/2} c_j \Re\left(\xi^{j - \frac{1}{2}}\right)
\end{align*}
being a purely imaginary quantity. Let us abbreviate
\[
  w_j := \Re\left(\xi^{j-\frac{1}{2}}\right) = \cos \left(\left(j-\tfrac{1}{2}\right)\tfrac{\pi}{n}\right).
\]
We can then write the problem of minimizing the constraint violation as
\[
  \min w^T c \quad \text{s.t. } w^T c \ge 0, c \in \{ \pm 1 \}^{n/2}.
\]
By substituting $c_j = 1 - 2 x_j$, we arrive at the SSP with real coefficients
\begin{equation} \label{eqn:SSP}
  \max w^T x \quad \text{s.t. } w^T x \le \frac{1}{2} \sum_{j=1}^{n/2} w_j, x \in \{0,1\}^{n/2}.
\end{equation}
SSPs are a special case of Knapsack Problems, which are known to be NP-hard. For an introduction, we refer the reader to~\cite{kellerer2004knapsack}. Fully polynomial-time approximation schemes are of no help here, because the coefficients $w$ are not rational, leading to extremely large constants that would lead to very extreme memory requirements.

\subsection{Phase~II: Code-restricted nonlinear problem}

Now we fix the code $c^\ast$ from Phase~I and solve MINLP~\eqref{eqn:minlp}, where it is computationally advantageous to reorder the sums in the constraints and to fix $\varphi_1$ and $\varphi_{n+1}$ to eliminate the continuous rotational symmetry. We obtain the NLP
\begin{equation} \label{eqn:nlp}
  \begin{alignedat}{3}
    \max~ & \sum_{j=1}^n 2 \sin \frac{\varphi_{j+1} - \varphi_{j}}{2} =: -f(\varphi) \\
    \text{s.t.~}
    & 0 = \sum_{j=2}^{n} (c^\ast_{j-1} - c^\ast_j) \cos \varphi_{j} - (c^\ast_1 + c^\ast_n) =: g_1(\varphi), \\
    & 0 = \sum_{j=2}^{n} (c^\ast_{j-1} - c^\ast_j) \sin \varphi_{j} =: g_2(\varphi), \\
    & 0 = \varphi_1 \le \varphi_2 \le \dotsb \le \varphi_{n+1} = \pi.
  \end{alignedat}
\end{equation}
Because the inequality constraints will not be active in the solutions we are interested in, we drop them in the remainder.
We consider a numerical method of Lagrange-Newton-type for solving NLP~\eqref{eqn:nlp}.

The Lagrangian corresponding to NLP~\eqref{eqn:nlp} without code symmetry exploitation reads
\[
  L(\varphi, y) = f(\varphi) + y_1 g_1(\varphi) + y_2 g_2(\varphi).
\]
For $j = 2, \dotsc, n$, its first derivatives are
\[
  \frac{\partial L}{\partial \varphi_j}(\varphi, y)
  = \frac{1}{2} \left(\cos \frac{\varphi_{j+1} - \varphi_{j}}{2}  - \cos \frac{\varphi_{j} - \varphi_{j-1}}{2} \right) + (c^\ast_{j-1} - c^\ast_{j}) (y_2 \cos \varphi_{j} - y_1 \sin \varphi_{j}).
\]
The Hessian matrix $H$ of $L$ with respect to $\varphi_2, \dotsc, \varphi_n$ is tridiagonal and symmetric with diagonal entries
\begin{align*}
  H_{jj}(\varphi, y)
  &= \frac{1}{4} \left( \sin \frac{\varphi_{j+1} - \varphi_{j}}{2} + \sin \frac{\varphi_{j} - \varphi_{j-1}}{2} \right) + (c^\ast_{j} - c^\ast_{j-1}) (y_1 \cos \varphi_{j} - y_2 \sin \varphi_{j})
\end{align*}
for $j = 2, \dotsc, n$ and off-diagonal entries
\[
  H_{j,j+1}(\varphi, y)
  = H_{j+1,j}(\varphi, y)
  = -\frac{1}{4} \sin \frac{\varphi_{j} - \varphi_{j+1}}{2}
\]
for $j = 2, \dotsc, n-1$.
The entries of the constraint Jacobian $C$ are for $j = 2, \dotsc, n$
\begin{align*}
  C_{1j}(\varphi, y) &= (c^\ast_{j} - c^\ast_{j-1}) \sin \varphi_{j}, &
  C_{2j}(\varphi, y) &= (c^\ast_{j-1} - c^\ast_{j}) \cos \varphi_{j}.
\end{align*}

Recalling that $\nabla_y L(\varphi, y) = g(\varphi)$, we can now perform a local Newton-type method on $w = (\varphi, y)$ according to
\[
  w^{k+1} = w^k - A_k \nabla L(w^k), \quad 
  w^{0} = 
  \begin{pmatrix}
    \varphi^\ast \\
    0
  \end{pmatrix}
  \quad \text{where } A_k^{-1} \approx K(w^k) :=
  \begin{pmatrix}
    H(w^k) & C(w^k)^T \\
    C(w^k) & 0
  \end{pmatrix}.
\]
It is important to use high-precision arithmetic to evaluate $\nabla L$. Inaccuracies in $A_k$ will reduce the convergence speed but will not change the fix points of the method, as long as $A_k$ stays invertible. Possible numerical realizations of $A_k$ are:
\begin{enumerate}
  \item Direct sparse decomposition: Evaluate $H(w^k)$ and $C(w^k)$, round to double precision, and use an off-the-shelf direct sparse decomposition algorithm in double precision on the rounded block matrix $K(w^k)$. This leads to a linearly convergent algorithm, where typically ten to fifteen digits of accuracy are gained per iteration.
  \item Schur complement: The symmetric tridiagonal $H(w^k)$ is easy to invert and the resulting Schur complement $-C(w^k) H^{-1}(w^k) C(w^k)^T$ is just a 2-by-2 matrix, which can also be easily inverted. All these operations can be performed in arbitrary precision arithmetic, which retains the locally quadratic convergence of the Newton method.
  \item MINRES: Use the Krylov method MINRES~\cite{paige1975solutions} in high-precision arithmetic as a direct solver by performing $n+1$ MINRES iterations, which also retains locally quadratic convergence.
  \item Simplified Newton variants: Take one of the first two approaches but always use the linearization in $w^0$, i.e., always use $A_0$ instead of $A_k$. This yields considerable speed-ups for each iteration. The rate of convergence, however, will be at most linear.
\end{enumerate}

We establish next that $K(w)$ is indeed invertible in a neighborhood of $w^0 = w^\ast$.

\begin{lemma} \label{lem:C_full_rank}
  If $0 < \varphi_2 < \dotsc < \varphi_n < \pi$, then $C(w)$ has rank two.
\end{lemma}
\begin{proof}
  Let us abbreviate the index set
  \[
    J = \{ j \in \{ 2, \dotsc, n \} \mid c^\ast_{j} \neq c^\ast_{j-1} \},
  \]
  which has at least two elements.
  A direct computation reveals
  \begin{align*}
    &\phantom{=}\,\det(C(w) C(w)^T)\\
    &= \sum_{i,j=2}^{n} (c^\ast_i - c^\ast_{i-1})^2 (c^\ast_j - c^\ast_{j-1})^2 \left( \sin^2 \varphi_i \cos^2 \varphi_j - \sin \varphi_i \cos \varphi_i \sin \varphi_j \cos \varphi_j \right) \\
    &= 16 \sum_{i,j \in J} \left( \sin^2 \varphi_i \cos^2 \varphi_j - \sin \varphi_i \cos \varphi_i \sin \varphi_j \cos \varphi_j \right) \\
    &= 16 \sum_{i \in J} \sum_{\substack{j \in J\\j > i}} \left( \sin^2 \varphi_i \cos^2 \varphi_j - 2 \sin \varphi_i \cos \varphi_i \sin \varphi_j \cos \varphi_j  + \sin^2 \varphi_j \cos^2 \varphi_i \right) \\
    &= 16 \sum_{i \in J} \sum_{\substack{j \in J\\j > i}} (\underbrace{\sin \varphi_i \cos \varphi_j - \sin \varphi_j \cos \varphi_i}_{=\sin (\varphi_i - \varphi_j) \neq 0} )^2 > 0.
  \end{align*}
  Hence, $C(w) C(w)^T$ has rank two, which is only possible if $C(w)$ has rank two.
\end{proof}
Lemma~\ref{lem:C_full_rank} also shows that any non-degenerate configuration of angles satisfies the Linear Independence Constraint Qualification for NLP~\eqref{eqn:nlp}.

\begin{lemma} \label{lem:H0_pos_def}
  The Hessian matrix $H(w^\ast)$ is positive definite for $(\varphi^\ast, 0)$.
\end{lemma}
\begin{proof}
  By the definition of $\varphi^\ast$, we have for all $j=2, \dotsc, n$ that
  \[
    \sin \frac{\varphi^\ast_j - \varphi^\ast_{j-1}}{2}
    = \sin \frac{\pi}{2n}.
  \]
  Hence,
  \[
    H(w^\ast)
    = \frac{\sin \frac{\pi}{2n}}{4}
    \begin{pmatrix}
      2 & -1 \\
      -1 & \ddots & \ddots \\
      & \ddots & \ddots & -1 \\
      & & -1 & 2
    \end{pmatrix},
  \]
  which has a complete set of eigenvectors $v^l \in \mathbb{R}^{n-1}$ of the form $v^l_j = \sin \frac{lj\pi}{n}$ because with the convention $v^l_0 = v^l_n = 0$, a trigonometric identity yields for all $l = 1, \dotsc, n-1$ that
  \[
    v^l_{j-1} + v^l_{j+1} = \sin \left(\frac{lj\pi}{n} - \frac{l\pi}{n} \right) + \sin \left(\frac{lj\pi}{n} + \frac{l\pi}{n} \right)
    = 2 v^l_j \cos \frac{l\pi}{n}.
  \]
  Hence, the eigenvalues of $H(w^\ast)$ are $\frac{1}{2} (1 - \cos \frac{l \pi}{n}) \sin \frac{\pi}{2n} > 0$ for $l = 1, \dotsc, n-1$ and $H(w^\ast)$ is positive definite.
\end{proof}
\begin{lemma} \label{lem:K_invertible}
  The block matrix $K(w)$ is invertible in a neighborhood of $w^\ast$.  
\end{lemma}
\begin{proof}
  The constraint Jacobian $C(w^\ast)$ has full rank and the Hessian $H(w^\ast)$ is positive definite by Lemmas~\ref{lem:C_full_rank} and~\ref{lem:H0_pos_def}. Hence $K(w^\ast)$ is invertible. 
  By continuity of $K$, it is invertible in a neighborhood of $w^\ast$.  
\end{proof}

Lemma~\ref{lem:K_invertible} assures us that the simplified Newton variants have well-defined iterates, even though convergence is only guaranteed if a solution of NLP~\eqref{eqn:nlp} is located sufficiently close to $w^\ast$.

The NLP~\eqref{eqn:nlp} can be reduced further by exploiting structure in the code $c^\ast$, e.g., code symmetry as in Mossinghoff's conjecture or by exploiting that pending diameters dissect the diameter graph cycle's angles equally. We shall see in Sec.~\ref{sec:numerics} that the real bottleneck of the two-phase approach is Phase~I, so that further improvements for Phase~II will not yield considerable computational improvements.

\section{Computational experience} \label{sec:numerics}

We want to provide information on the computational experiences we gained with the approaches put forward in Sec.~\ref{sec:two_phase}.

\subsection{Phase~I}

Up to $n=32$, we enumerated all possible codes (without the additional symmetry of Mossinghoff's conjecture) using the C program accompanying the paper~\cite{sawada2001generating}.

We also experimented with minimizing the constraint violation in MINLP~\eqref{eqn:minlp} for fixed $\varphi=\varphi^\ast$ using the angular constraints described at the beginning of Sec.~\ref{sec:phase_one_power_of_two} in Gurobi~\cite{gurobi}, which works well up to $n=32$.

For $n \ge 64$, the default settings of Gurobi lead to acceptance of non-optimal iterates, which is due to double precision arithmetic incurring too large numerical errors. Experiments with tweaking the options \texttt{NumericFocus} to its maximum of 3 and enabling \texttt{Quad} for quadruple precision in the Simplex method were not enough to obtain reliable results. This applies also for the SSP problem~\eqref{eqn:SSP} resulting from assuming Mossinghoff's conjecture to be true.

Using the exact rational version of SCIP~\cite{bestuzheva2023enabling,cook2013hybrid} for extra precision to solve the SSP~\eqref{eqn:SSP} was also unsuccessful. With default parameter settings, the solver ran into memory issues on a large machine with 500\,GB of RAM. These memory issues can be mitigated by deactivating cut generation (option \texttt{separating/maxcuts=0}) and changing the node selection to depth-first-search (option \texttt{nodeselection/dfs/stdpriority=1000000}). The performance is, however, quite slow due to the usage of rational arithmetic and it appears that a large part of the branching tree must be visited.

Hence, we implemented two concise C programs, which can be obtained from the authors on reasonable request. The first uses GCC's \texttt{libquadmath} quadruple floating point arithmetic, the second uses purely 128-bit integer arithmetic for rational numbers with a fixed common denominator of $2^{128} / n$. The chosen denominator is sufficiently large because the sum of all $w_j$ is bounded above by $n/2$. This choice yields an absolute accuracy of $3.8 \cdot 10^{-37}$. Both programs implement a basic recursive enumeration of all possible solutions of SSP~\eqref{eqn:SSP} with early branch pruning if the currently best objective cannot be reached or if the constraint is already violated. The programs also provide a simple means of parallelization by fixing a code suffix of given length. For instance, a parallelization with four processors can then be obtained by calling the program four times with the suffixes \texttt{00}, \texttt{01}, \texttt{10}, \texttt{11}.

On one core of an AMD EPYC 7543 32-Core 3.3 GHz CPU, the quadruple precision version takes 9.3 seconds for the complete enumeration for $n=64$, while the 128-bit integer version takes only 0.95 seconds. The enumeration programs clearly beat the LP-based approaches we tried.

We present the resulting gaps $\frac{1}{2}\sum_{j=1}^{n/2} w_j - w^T x$ of the constraint of SSP~\eqref{eqn:SSP}, the runtime of the 128-bit integer program, and the corresponding quarter code in Tab.~\ref{tab:codes}. Due to the combinatorial explosion of the number of codes, we were not able to finish the Phase~I computations for $n=128$. However, we were able to find a code with very small gap within two days of parallel computations. A back-of-the-envelope calculation leads to an estimated runtime of a staggering 14 years for the 128-bit integer variant of the Phase~I program on 64 cores in parallel.

\begin{table}[htbp]
  \centering
  \begin{tabular}{cccc}
    $n$ & Gap & Runtime [s] & Quarter code $c_1, \dotsc, c_{n/2}$ \\
    \hline
    4   & $6.533 \cdot 10^{-1}$  & 0.00 & \texttt{+-} \\
    8   & $3.007 \cdot 10^{-1}$  & 0.00 & \texttt{+--+} \\
    16  & $2.070 \cdot 10^{-2}$  & 0.00 & \texttt{+--+-++-} \\
    32  & $3.409 \cdot 10^{-5}$  & 0.00 & \texttt{+-++--+-+-+---++} \\
    64  & $1.984 \cdot 10^{-9}$  & 0.95 & \texttt{-++++++-----+--+-+++--+---+--+++}\\
    128 & $2.000 \cdot 10^{-16}$ & $-^\ast$ & \texttt{-+++---+-++++-+++++----+-+---+-+} \\
    &&& $\hookrightarrow$ \texttt{+----++++-+-----+--+-+-+--+--+--} \\
    \hline
  \end{tabular}
  \caption{Quarter codes computed with the two-phase approach. $^\ast$The code used for $n=128$ is the best one found after two days of running the Phase~I program on 64 cores and is probably suboptimal.}
  \label{tab:codes}
\end{table}

\subsection{Phase~II}

We implemented the Lagrange-Newton method for NLP~\eqref{eqn:nlp} in a Python notebook using the GNU Multiple Precision Arithmetic Library via the Python package \texttt{gmpy2} and sparse direct decomposition of $K(w^k)$ in double precision using UMFPACK~\cite{davis2004umfpack} via \texttt{scikits-umfpack}. We ran the Lagrange-Newton method on codes computed in Phase~I. For $n \ge 32$, we used the simplified Newton variant, in which $A_k = A_0$ in every iteration. The number of required iterations and runtimes using an arithmetic precision of 360 binary digits and a Newton tolerance of 320 binary digits ($\approx 4\cdot 10^{-103}$) on the Newton increment norm are displayed in Tab.~\ref{tab:phase_two}. The runtimes were obtained by running the (simplified) Newton method five times and taking the minimum of the measured times on a 2 GHz Quad-Core Intel Core i5 MacBook Pro. The resulting polygons for $n \ge 32$ are depicted in Fig.~\ref{fig:polys_vary_n} and the corresponding objectives are given in Tab.~\ref{tab:perimeters_vary_n}. We can see that for large $n$, the simplified Newton method performs highly efficiently.

The runtimes also show that the real bottleneck of the two-phase approach is Phase~I. Additional structure exploitation, e.g., through reduced versions of NLP~\eqref{eqn:nlp} is not necessary.

\begin{table}[htbp]
  \centering
  \begin{tabular}{ccc}
    $n$ & (Simplified$^\ast$) Newton it. & Runtime [s] \\
    \hline
    4   & 10 & 0.02 \\
    8   & 10 & 0.03 \\
    16  &  9 & 0.16 \\
    32  & $17^\ast$ & 0.09 \\
    64  & $9^\ast$ & 0.02 \\
    128 & $8^\ast$ & 0.05 \\
    \hline
  \end{tabular}
  \caption{Number of Phase~II Newton iterations and corresponding runtime for different codes of length $n$ from Phase~I. For $n \ge 32$, a simplified Newton method was used, indicated by $^\ast$.}
 \label{tab:phase_two}
\end{table}




\section{The longest perimeter of a small octagon is algebraic} \label{sec:perim_algebraic}

Phase~I delivers a code $c^\ast$. As an alternative to Phase~II, it is in some cases possible to derive the solution of the maximum perimeter problem for the fixed code $c^\ast$ by analytical means.
 
\subsection{Statements of the results}

%
%
%
%
As mentioned in the introduction, the maximum area of a small hexagon received quite some attention.
It occurs as \emph{Graham's hexagon constant} in the fundamental encyclopedia \cite[p.~526]{Finch2003}
and is available in Mathematica as \texttt{n=6; GraphData[{"BiggestLittlePolygon", n}]}.%
\footnote{But be aware, the octagon data for $n=8$ are quite inaccurate.}
To some extent this popularity is based on its representation as an algebraic number.%
\footnote{It seems to be worthwhile to cite the comment by H. Hadwiger in \cite{bieri1961ungeloeste}: Nun ist es H. Bieri (Bern) gelungen, das Problem für $n=6$ im achsensymmetrischen Fall zu lösen. Hierbei ist einzuräumen, dass das extremale Polygon auch im allgemeinen Fall achsiale Symmetrie aufweisen dürfte, so dass die gefundene Lösung vermutlich allgemeine Gültigkeit hat. Die von Bieri erzielte Lösung ist indessen nicht einfach angebbar, und die Form, in der die eindeutige Kennzeichnung des extremalen Polygons und des maximalen Flächeninhaltes erfolgen muss, regt an, sich über den Begriff Lösung eines elementargeometrischen Problems einige nützliche Gedanken zu machen.}
This suggests asking for similar results for the longest perimeter of a small octagon. 

Here we demonstrate that the square of the longest perimeter of 
a small octagon is a root of an integer polynomial of degree $48$.
We obtained this polynomial by three independent approaches. The
derivations are based on a formulation using cartesian coordinates or via
a trigonometric formulation with half angle substitution.
The most efficient approach uses \name{Lagrange} multipliers.
The application of this idea to \name{Graham}'s hexagon is somewhat
hidden as an example in \cite[p.~46]{trott_mathematica_2006}.
Therefore, we begin with revisiting the small hexagon and
the small octagon of largest area.
We also consider the case of equilateral small octagons and derive a polynomial of degree $6$ for the corresponding longest perimeter.
This polynomial has already been mentioned without proof in \cite{Bingane2022a}. 

Our derivations rely on algorithms dealing with
ideals of multivariate polynomials with integer coefficients. We carried
out all computations using implementations of such algorithms 
in \name{SageMath}
\cite{Sage} as well as in \name{Mathematica} \cite{Mathematica}.
We only provide brief descriptions of our approach, which are not to be considered as complete
proofs.

\begin{proposition}
    The square of the maximum perimeter $p_8 \approx 3.12114713405983$ of a 
    small octagon is a root of an integer polynomial $P_8$ of degree $48$,
    see \autoref{sec:app}.
    The maximum perimeter $e_8= 8 s_8 \approx  3.095609317476962$ of an equilateral 
    small octagon belongs to the side length $s_8$, which is a root of the integer polynomial 
    \begin{align*}
        E_8(t)&=2t^6 - 18t^5 + 57t^4 - 78t^3 + 46t^2 - 12t + 1
    \end{align*}
  of degree $6$.
  \end{proposition}

We include a formulation of the results concerning the maximum area, since we are going to
present alternative derivations. 

\begin{proposition}[\cite{graham_largest_1975,audet2020using}]
    The maximum area $a_6\approx 0.6749814429301$ of a 
    small hexagon is a root of the integer polynomial
    \begin{align*}
        A_6(t)&=4096 t^{10}+8192 t^9-3008 t^8-30848 t^7+21056 t^6+146496 t^5\\
        &\qquad\qquad -221360t^4+1232t^3+144464 t^2-78488 t+11993    
    \end{align*}
    of degree $10$.
    The maximum area $a_8\approx 0.726868482751$ of a 
    small octagon is a root of an
    integer polynomial $A_8$ of degree $42$, see \cite{audet2020using}.
\end{proposition}

%

\subsection{The largest area small hexagon and octagon}

\subsubsection{\name{Lagrange} multipliers for the nonsymmetric hexagon}

\def\bmscale{0.85}
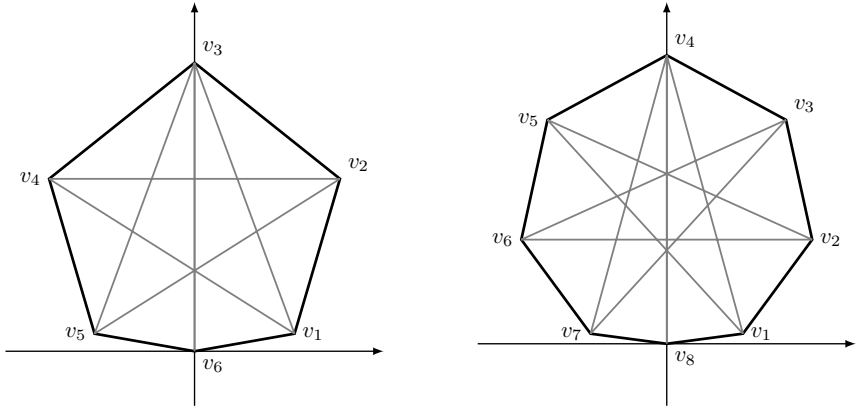
\begin{figure}
	\hspace*{-25mm}
    {     \scalebox{\bmscale}{




\definecolor{TUCgreen}{rgb}{0,0.55,0.31}
\definecolor{TUCgrey1}{rgb}{0.5,0.5,0.5}                                                   
\definecolor{TUCgrey2}{rgb}{0.9,0.9,0.9}
\definecolor{TUCred}{rgb}{0.55,0.11,0}
\definecolor{green}{rgb}{0,0.55,0.31}
\definecolor{grey}{rgb}{0.9,0.9,0.9}
\definecolor{red}{rgb}{0.55,0.11,0}

\newcommand{\tscale}{4.5}
\newcommand{\vt}{\vartheta}
\newcommand{\vp}{\varphi}
\newcommand{\ve}{\varepsilon}
\newcommand{\twod}{\mathinner{\ldotp\ldotp}}
\newcommand{\rg}[2]{\{#1\mathbin{\upupharpoons}#2\}}
\newcommand{\cj}[1]{\overline{#1}}


\def\lw{0.6}
\def\xmax{0.6}
\def\ymax{1.1}
\def\ymin{-0.2}
\def\R{1}
\def\edto(#1){edge ++(#1) to ++(#1)}
\tikzstyle{vecr}=[->,very thick,red,line cap=round]
\tikzstyle{vecg}=[->,very thick,green]
\tikzstyle{vecrd}=[->,very thick,red,dashed,line cap=round]
\tikzstyle{vecgd}=[->,very thick,green,dashed,line cap=round]
\tikzstyle{vecb}=[->,very thick,blue,line cap=round]
\tikzstyle{lir}=[-,thick,white!50!black]
\tikzstyle{lib}=[-,very thick,black]

  \def\vtclasp{(0.262, 0.),(0.5, 0.326),(0.41, 0.741),(0., 0.965),(-0.41, 
  0.741),(-0.5, 0.326),(-0.262, 0.),(0., -0.035)}

  \def\vtchex{(0.343771, 0.), (0.5, 0.536703), (0., 0.939053), (-0.5, 
  0.536703), (-0.343771, 0.), (0., -0.0609467)}

 \begin{tikzpicture}[scale=4.5]
    \foreach \coord [count=\ind] in \vtchex{
         \coordinate [
         at=\coord] (v_\ind);}
    \coordinate (O) at (0,0); 
  \draw[->,line width=\lw] ($(-\xmax-0.05,0)+(v_6)$) -- ($(\xmax+0.05,0)+(v_6)$);
  \draw[->,line width=\lw] (0,\ymin-0.05) -- (0,\ymax+0.05);
    \draw[lib] plot[samples at={1,...,6}] (v_\x) -- cycle;
    \draw[lir] (v_1)--(v_3) (v_2)--(v_4) (v_3)--(v_6)
    (v_3)--(v_5) (v_1)--(v_4) (v_2)--(v_5);
    \foreach \ind in {2,...,3}
      \draw (v_\ind) node[above right] {$v_\ind$};
    \foreach \ind in {4,...,5}
      \draw (v_\ind) node[left] {$v_\ind$};
    \draw (v_1) node[right] {$v_1$};
    \draw (v_6) node[below right] {$v_6$};
%
\end{tikzpicture}

}
  \hspace*{-20mm}
		\scalebox{\bmscale}{




\definecolor{TUCgreen}{rgb}{0,0.55,0.31}
\definecolor{TUCgrey1}{rgb}{0.5,0.5,0.5}                                                   
\definecolor{TUCgrey2}{rgb}{0.9,0.9,0.9}
\definecolor{TUCred}{rgb}{0.55,0.11,0}
\definecolor{green}{rgb}{0,0.55,0.31}
\definecolor{grey}{rgb}{0.9,0.9,0.9}
\definecolor{red}{rgb}{0.55,0.11,0}

\newcommand{\tscale}{4.5}
\newcommand{\vt}{\vartheta}
\newcommand{\vp}{\varphi}
\newcommand{\ve}{\varepsilon}
\newcommand{\twod}{\mathinner{\ldotp\ldotp}}
\newcommand{\rg}[2]{\{#1\mathbin{\upupharpoons}#2\}}
\newcommand{\cj}[1]{\overline{#1}}


\def\lw{0.6}
\def\xmax{0.6}
\def\ymax{1.1}
\def\ymin{-0.2}
\def\R{1}
\def\edto(#1){edge ++(#1) to ++(#1)}
\tikzstyle{vecr}=[->,very thick,red,line cap=round]
\tikzstyle{vecg}=[->,very thick,green]
\tikzstyle{vecrd}=[->,very thick,red,dashed,line cap=round]
\tikzstyle{vecgd}=[->,very thick,green,dashed,line cap=round]
\tikzstyle{vecb}=[->,very thick,blue,line cap=round]
\tikzstyle{lir}=[-,thick,white!50!black]
\tikzstyle{lib}=[-,very thick,black]

  \def\vtclasp{(0.262, 0.),(0.5, 0.326),(0.41, 0.741),(0., 0.965),(-0.41, 
  0.741),(-0.5, 0.326),(-0.262, 0.),(0., -0.035)}

  \def\vtchex{(0.343771, 0.), (0.5, 0.536703), (0., 0.939053), (-0.5, 
  0.536703), (-0.343771, 0.), (0., -0.0609467)}

 \begin{tikzpicture}[scale=4.5]
    \foreach \coord [count=\ind] in \vtclasp{
         \coordinate [
         at=\coord] (v_\ind);}
    \coordinate (O) at (0,0); 
  \draw[->,line width=\lw] ($(-\xmax-0.05,0)+(v_8)$) -- ($(\xmax+0.05,0)+(v_8)$);
  \draw[->,line width=\lw] (0,\ymin-0.05) -- (0,\ymax+0.05);
    \draw[lib] plot[samples at={1,...,8}] (v_\x) -- cycle;
    \draw[lir] (v_1)--(v_5) (v_2)--(v_6) (v_3)--(v_7) (v_4)--(v_8)
    (v_3)--(v_6) (v_2)--(v_5) (v_1)--(v_4) (v_7)--(v_4);
    \foreach \ind in {1,...,2}
      \draw (v_\ind) node[right] {$v_\ind$};
    \foreach \ind in {3,...,4}
      \draw (v_\ind) node[above right] {$v_\ind$};
    \foreach \ind in {5,...,7}
      \draw (v_\ind) node[left] {$v_\ind$};
    \draw (v_8) node[below right] {$v_8$};
%
\end{tikzpicture}

} }
\caption{The small hexagon and the small octagon of largest area. Both are symmetric through the
(vertical) pending diameter chord. \label{fig:area}}
\end{figure}

The enumeration of the vertices is illustrated in \autoref{fig:area}.

The vertices of the hexagon are parametrized by
\begin{align*}
    v_1&=(x_1,y_1),\quad v_2=(x_2,y_2), \quad v_3=(0,1), \\
     v_4&=(x_4,y_4), \quad v_5=(-x_1,y_1), \quad v_6=(0,0).
\end{align*} 
The symmetry of $v_1$ and $v_5$ incorporates the fact that the pending diameter chord
$v_3v_6$, which we fix, needs to be perpendicular to $v_1v_5$ for maximum area.
We also emphasize that we do not require full symmetry through the line $v_6v_3$. In fact, our derivation
shows this symmetry property of the extremal hexagon to hold.
Expressing the area $t$ in terms of the variables gives the problem
\begin{align*}
    &\max_{x_1,y_1,x_2,y_2,x_4,y_4} t=\big(x_2 + x_4 (y_1-1) - x_2 y_1 + x_1 (y_2 + y_4)\big)/2 \quad \text{subject to}\\
    &\qquad\qquad \abs{v_1-v_3}^2=\abs{v_1-v_4}^2=\abs{v_2-v_4}^2=\abs{v_2-v_5}^2=1. 
    \end{align*}

We obtain the \name{Lagrange} function
\begin{align*}
L&= t +\lambda_1 (\abs{v_1-v_3}^2-1)+\lambda_2 (\abs{v_1-v_4}^2-1)\\
&\qquad\qquad\qquad\qquad +\lambda_3 (\abs{v_2-v_4}^2-1)+\lambda_4 (\abs{v_2-v_5}^2-1)
\end{align*} 
and the system
\begin{align}
	\frac{\partial L}{\partial x_1}=0,\; \dots,\; 
	\frac{\partial L}{\partial \lambda_4}=0.
\end{align}
This resulting \name{Lagrange} system is a system of polynomial equations in the variables
$t,x_1,y_1,x_2,y_2,x_4,y_4$ and four multipliers.
The basic polynomial in the single variable $t$, that means in the area, is computed by
\texttt{GroebnerBasis} as
\begin{align}
	B(t) = A_6(-t) A_6(t) Q(t),
\end{align}
where $Q$ is of degree $18$ in $t^2$.

The symmetry assumption would reduce the problem by explicitly setting $x_2=-x_4=\tfrac{1}{2}$, $y_4=y_2$, and
$\lambda_3=\lambda_4=0$, and would give the polynomial $A_6$ without any factorization. 

%

\subsubsection{\name{Lagrange} multipliers for the symmetric octagon}
Here we completely exploit the known symmetry and incorporate it in the problem formulation.
The vertices of the octagon are
\begin{align*}
  v_1 &= (x_1, y_1), \quad v_2 = (\tfrac{1}{2}, y_2),\quad  v_3 = (x_3, y_3),\quad  v_4 = (0, 1),\\
 v_5 &= (-x_3, y_3),\quad  v_6 = (-\tfrac{1}{2}, y_3),\quad  v_7 = (-x_1, y_1),\quad  v_{8}=(0,0),
\end{align*}
its area is then given by
\begin{align*} 
t=x_3 - x_3 y_2 + (-y_1 + 2 x_1 y_2 + y_3)/2
\end{align*}
and there remain the constraints 
$\abs{v_1-v_4}^2=\abs{v_1-v_5}^2=\abs{v_2-v_5}^2=\abs{v_2-v_6}^2=1$.
The corresponding Lagrange system is then reduced to the polynomial, which was already derived
in \cite{audet2020using}.

%
%
%

\subsection{The longest perimeter small octagon}

\def\bmscale{0.85}
\begin{figure}
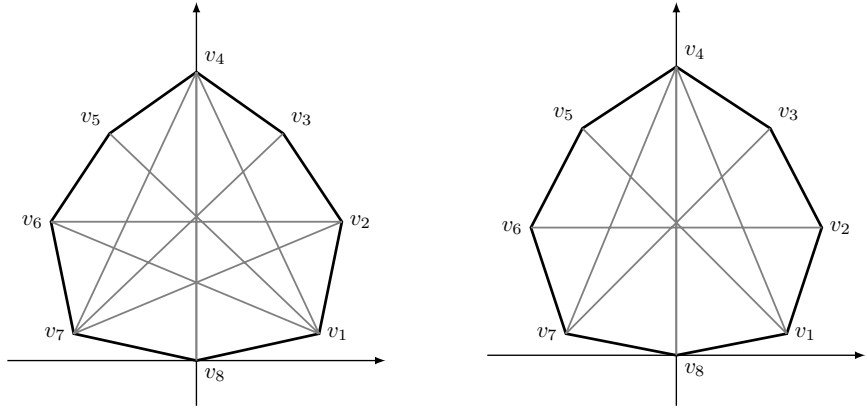

  \centering
  \hspace*{-30mm}
    {     \scalebox{\bmscale}{\input{lpsp_octagon.tex}}
    \hspace*{-20mm}
		\scalebox{\bmscale}{\input{equilateral_lps_octagon}} }
\caption{The small octagon and the equilateral small octagon of longest perimeter. Both are symmetric, the equilateral solution has only six diameter chords. \label{fig:perimeter}}
\end{figure}

As evident from \cite{audet2007small} and from our computations, the extremal octagon $P$ has 
eight diagonals of unit length. This diameter graph consists
of a cycle of five edges and three pending edges. The octagon $P$ is symmetric
through one of the pending edges. For an illustration and our notation
see \autoref{fig:perimeter}.

\subsubsection{Half-angle substitution}
We consider the angles $\alpha=\frac12 \angle{v_4v_1v_6}$ and
$\beta=\frac14 \angle{v_1v_4v_7}$, which implies $\frac12 \angle{v_6v_2v_7}=\pi/4-\alpha-\beta$.
Now, the extremal problem is formulated as
\begin{alignat*}{3}
    \max_{\alpha,\beta}~ & 8\sin\frac\alpha2+4\sin\beta
    +4\sin(\pi/4-\alpha-\beta) \\
    \text{s.t.~}
    &2\sin(2\beta)+2(1-\cos(\pi/2-2\alpha-2\beta)) = 1.
\end{alignat*}
Note that the constraint is related to the length $|v_6-v_2|=1$.
Using the substitutions
$x = \tan(\alpha/4)$, $y=\tan(\beta/2)$, this problem becomes a rational system.
The \name{Mathematica}
functions \texttt{Resultant} and \texttt{Discriminant} deliver
a polynomial of degree 48 in $t$ in integers extended by
$\sqrt{2}$, which is equivalent to the one given in Appendix~\ref{sec:app}.

\subsubsection{\name{Lagrange} multipliers}
Respecting the symmetry through the pending diameter chord $v_8v_4$,
the vertices of $P$ are shifted vertically downwards and parameterized by
\begin{align*}
    v_1&=(x,0),\quad v_2=(\tfrac{1}{2},v_2),\quad v_4=(0,1+y_8), \\
    v_6&=(-\tfrac{1}{2},y_2), \quad v_7=(-x,0), \quad v_8=(0,y_8).
\end{align*}
The vertex $v_3$ has to be the midpoint of the arc $v_2v_4$
of the circle around $v_7$ of unit radius. Introducing the 
slack variables $h_1=\abs{v_1-v_8}$, $h_2=\abs{v_2-v_1}$,
$h_3=\abs{v_3-v_2}=\abs{v_4-v_3}$, this implies the relation
$\abs{v_4-v_2}^2+h_3^4-4h_3^2=0$. Hence, the extremal problem
may be written as
\begin{alignat*}{3}
    \max_{x,y,v}~& 2(h_1+h_2+2h_3)\\
    \text{s.t.~} &
    \abs{v_1-v_4}^2=\abs{v_2-v_7}^2=1,\quad
    \abs{v_1-v_8}^2=h_1^2, \quad \abs{v_2-v_1}^2=h_2^2, \\
    & \abs{v_4-v_2}^2+h_3^4-4h_3^2=0.
\end{alignat*}
The \name{Lagrange} system is a polynomial system in the variables
$p,x,y_2,y_8,h_1,h_2,h_3$ and five multipliers.
The function \texttt{GroebnerBasis} in \name{Mathematica}
easily comes up with the announced polynomial for the perimeter $t$.

As described above, we elaborated on a purely numerical two-phase approach for computing the maximum perimeters of small polygons. 
As a sanity check, the results
for the root $p_8$ of the integer polynomial are identical with the ones reported in Tab.~\ref{tab:perimeters_n_8}.

\section{Conclusions} \label{sec:conclusions}

We derived two new MINLP formulations for computing lower bounds for Reinhardt's perimeter problem, which is still an open problem after more than 100 years. We argue that the derived lower bounds are strict if it is true that optimal polygons have a complete set of diameters. We propose a two-phase approach that consists of a (generalized) SSP to compute the best code, for which the regular $n$-gon, which realizes Reinhardt's upper bound, is almost feasible. Then we refine the regular $n$-gon in the second phase to a feasible polygon with large perimeter using an arbitrary precision Lagrange-Newton-type method. We share ample computational experience and provide high-accuracy solution candidates up to $n=64$ and a probably suboptimal candidate for $n=128$ with large perimeter. We further show that the perimeter of the octagon with largest perimeter is an algebraic number and provide the corresponding integer polynomials of degree 48.

\backmatter


\section*{Declarations}

The authors declare that they have no known competing financial or non-financial interests that could have appeared to influence the work in this paper.

\begin{appendices}
  
  
  \section{The polynomial of the longest perimeter small octagon}\label{sec:app}
  Here we provide the integer coefficients of the polynomial derived in Sec.~\ref{sec:perim_algebraic}, which can also be used to computate the longest perimeter as a root of this polynomial to arbitrary precision.

 \everydisplay{\fontsize{8pt}{10pt}\selectfont}
        
\allowdisplaybreaks
   \begin{align*}
    q_8(\, t)=&\, t^{48}-2688 \, t^{47}+3446272 \, t^{46}-2807592960 \, t^{45}
    +1633674266624 \, t^{44}\\ 
    &-723588467449856 \, t^{43}+2538898109167370,24 \, t^{42}\\ &-72517854035614629888 \, t^{41}+17195167136435692371968 \, t^{40}\\ &-3434820210452081718329344 \, t^{39}+584574334125421235902873600 \, t^{38}\\ &-85517739615251895256392663040 \, t^{37}
    +10829446302748698450583070179328 \, t^{36}\\ &-1193816644613689654665612818907136 \, t^{35} +115086843866723348026033261837811712 \, t^{34} \\ &-9738082617441544437635181820812197888 \, t^{33}\\ &+725423774155626603530948947333269684224 \, t^{32} \\ &-47694399655201711781649368889348291821568 \, t^{31}\\ &+2773474766447168322184997293630827538677760 \, t^{30}\\ &-142915937585753382520287017253466824236859392 \, t^{29}\\ &+6537352911881897599073318260145551517258088448 \, t^{28}\\ &-265914881136849371015171799699357177396926087168 \, t^{27}\\ &+9635371582576344208603139477144787322898164482048 \, t^{26}\\ &-311558095037500866248196651997764943419255015079936 \, t^{25}\\ &+9003894302377052774819990635603424836273155562536960 \, t^{24}\\ &-232798593505589953724659815997148521459619251942850560 \, t^{23}\\ &+5384698849214794614511096063213254330174249732122083328 \, t^{22}\\ &-111229883050697265095697776791662655429798099484617474048 \, t^{21}\\ &+2043943451881331447819048584603121094991258075349364244480 \, t^{20}\\ &-33204170232890601212009526399439211418502869355381160673280 \, t^{19}\\ &+472940328670015419031704476959013624500278043685261018136576 \, t^{18}\\ &-5852634224113275457409744736043691759211050964981756612050944 \, t^{17}\\ &+62479370348335228971031421065092236129622542636279079783890944 \, t^{16}\\ &-576375293461315675942032249195563488574543633050881901599916032 \, t^{15}\\ &+4710488377079952771836133728493511544552402556895090911062523904 \, t^{14}\\ &-36368623222587641069111818667387388732204823767000547849301655552 \, t^{13}\\ &+286657343740444279969083463862811708030275860108457118864141975552 \, t^{12}\\ &-2300409752638978585592882084275762121478341006958876960731097464832 \, t^{11}\\ &+17022846567049524499840649827179007136206891115523021868856676188160 \, t^{10}\\ &-105271885062338881752377059792380965744985216725352909959696593453056 \, t^9\\ &+519661898036620827612298078628038297122538078160606821748651700781056 \, t^8\\ &-2058218618695242650197271615605589498521465900890487961554567777222656 \, t^7\\ &+6790786688361483921080280822642541944806188768497117715297100944637952 \, t^6\\ &-18953181950142721458054232615448184781228851292920395499757235462995968 \, t^5\\ &+42033580403345828191224614534012438972659544181809655477758966164881408 \, t^4\\ &-66202009827628255733947076652696118184733299289788671027802646213820416 \, t^3\\ &+66152337375333098717233438245411622983399457619149468736826967974739968 \, t^2\\ &-36757562971957762256433166768884793509267603912864060192884364027101184 \, t\\ &+8563803984117479982232573393172818348276610047246959316838387256655872
   \end{align*}


\end{appendices} 


\end{document}